\documentclass{article}
\usepackage[utf8]{inputenc}
\usepackage{amsthm,amsmath,amssymb,xcolor}
\usepackage{tikz,tkz-berge}
\usepackage{multicol}
\usepackage{enumerate}
\usetikzlibrary{arrows}
\usepackage{hyperref}
\usepackage{tabularx}
\usepackage[export]{adjustbox}
\usepackage{appendix}
\usepackage[normalem]{ulem}
\usepackage{authblk}

\newtheorem{thm}{Theorem}

\newtheorem{cor}[thm]{Corollary}

\newtheorem{lemma}[thm]{Lemma}
\newtheorem{prop}[thm]{Proposition}
\theoremstyle{definition}
\newtheorem{Ex}{Example}

\renewcommand\det{\operatorname{Det}}
\newcommand\dist{\operatorname{Dist}}
\newcommand\aut{\operatorname{Aut}}

\definecolor{lblue}{RGB}{186,225,255}
\definecolor{lblue2}{RGB}{100,175,255}

\definecolor{lblue}{RGB}{186,225,255}
\definecolor{lblue2}{RGB}{100,175,255}
\usepackage{todonotes}
\definecolor{spurple}{RGB}{221,160,221}
\definecolor{pgreen}{RGB}{153,255,204}
\definecolor{korange}{RGB}{255,204,102}

\definecolor{spurple}{RGB}{221,160,221}
\definecolor{pgreen}{RGB}{153,255,204}
\definecolor{korange}{RGB}{255,204,102}

\newcommand\Z{{\mathbb Z}}

\begin{document}

\title{Symmetry Parameters of Two-Generator Circulant Graphs}

\author[1]{Sally Cockburn} 
\author[2]{Sarah Loeb}
\affil[1]{\url{scockbur@hamilton.edu}, Hamilton College, Clinton, NY}
\affil[2]{\url{sloeb@hsc.edu}, Hampden-Sydney College, Hampden-Sydney, VA}

\renewcommand\Affilfont{\footnotesize}

\date{\today}

\maketitle

\begin{abstract}
 The derived graph of a voltage graph consisting of a single vertex and two loops of different voltages is a circulant graph with two generators. We characterize the automorphism groups of connected, two-generator circulant graphs, and give their determining and distinguishing number, and when relevant, their cost of 2-distinguishing. We do the same for the subdivisions of connected, two-generator circulant graphs obtained by 
 replacing one loop in the voltage graph with a directed cycle.
\end{abstract}



\section{Introduction}

\begin{figure}[h]
\includegraphics[width=0.6\textwidth, center]{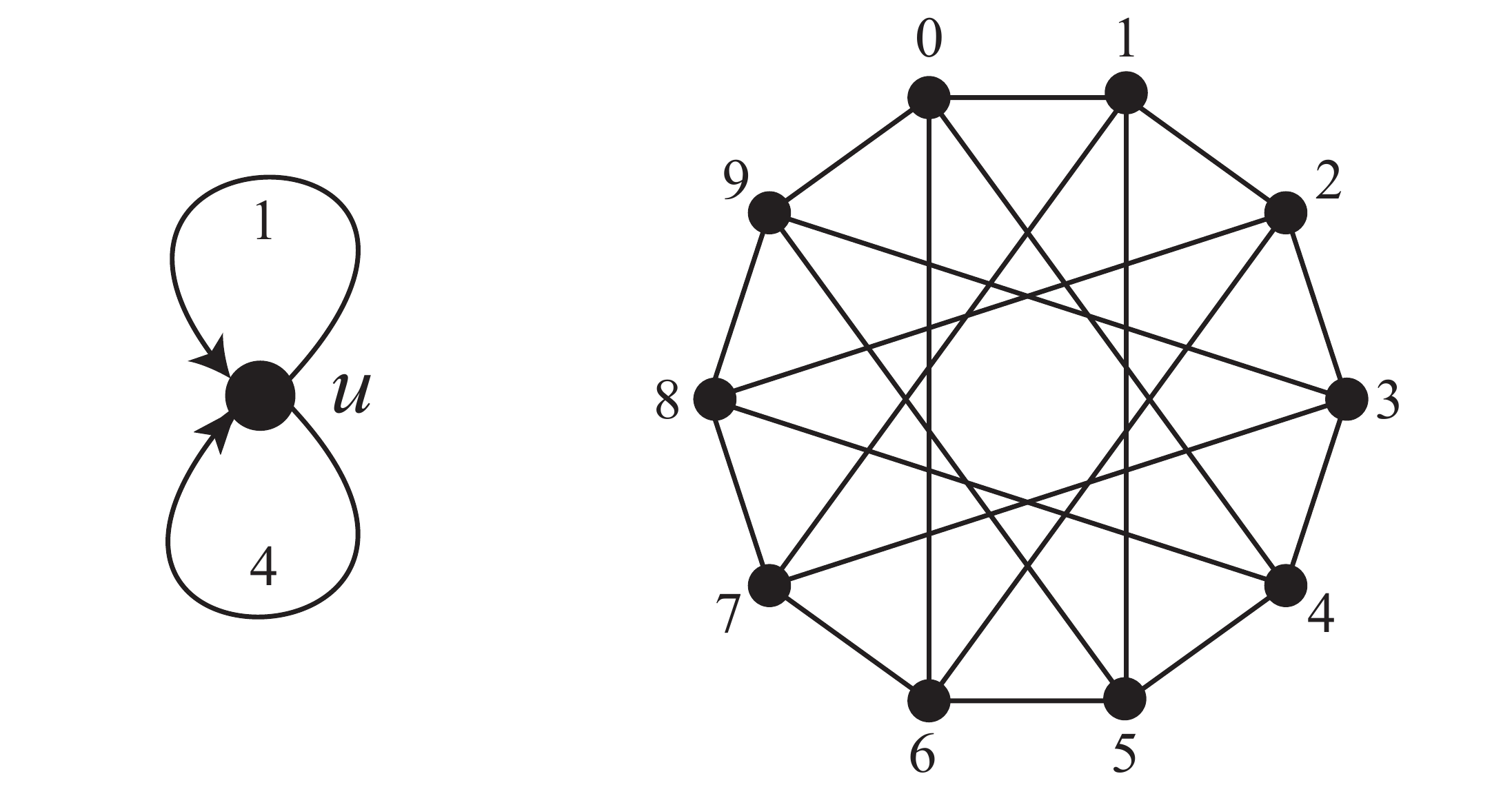}
\caption{Voltage graph $(B_2, \Z_{10})$ and derived graph $C_{10}(1, 4)$.}\label{fig:EZegVoltage}
\end{figure}

A voltage graph consists of a base directed graph $D= (V, E)$, a group $\Gamma$, and a voltage function $\phi:E \to \Gamma$. The associated derived directed graph $D^\phi$ has vertex set $\{u_a \mid u \in V, a \in \Gamma\}$ and arc set $\{e_a \mid e \in E, a \in \Gamma\}$; if $e = (u,v)$ and $\phi(e) = b$, then $e_a=(u_a, v_{ab})$. For more background on the voltage graph construction, see~\cite{GT1987}. 
 A particularly simple example has 
 a base directed graph consisting of a single vertex $u$ and two directed loops (sometimes called a bouquet), denoted $B_2$, and group $\Gamma =\Z_n$. Because the base graph has only one vertex, we can denote vertices in the derived graph simply as elements of $\Z_n$. We will denote an element of $\Z_n$ with an integer representative; for $a, b \in \Z$, we use the notation $a \equiv b$ to denote equality of the equivalence classes in $\Z_n$.
 The voltages on the two loops are denoted $i$ and $j$. Then the underlying undirected graph of the associated derived graph has vertex set $\Z_n$, with $a, b \in \Z_n$ adjacent if and only if $a-b \equiv \pm i
 $ or $a-b \equiv \pm j$. This is the circulant graph with two generators, commonly denoted $C_n(i,j)$. Figure~\ref{fig:EZegVoltage} shows an example with $n=10$, where the two loops have voltages $i=1$ and $j = 4$. 

We present some facts about two-generator circulant graphs. It is clear from the definition that $C_n(i,j) = C_n(n-i, j) = C_n(i, n-j) = C_n(n-i, n-j)$, so throughout this paper, we will assume that $0 < i < j \le n/2$. 
If $k$ is a unit in $\Z_n$, then multiplying vertices by $k$ is a graph isomorphism $C_n(i, j) \cong C_n(ik, jk)$. In particular, if $i$ is a unit in $\Z_n$, then $C_n(i,j)\cong C_n(1, i^{-1}j)$ (and if $j$ is a unit, then $C_n(i,j)\cong C_n(1, ij^{-1})$). In what follows, we will assume that either $i=1$ or neither $i$ nor $j$ is a unit in $\mathbb Z_n.$
As noted in~\cite{Brooks2021}, $C_n(i,j)$ is connected if and only if $\gcd(n,i,j) = 1$. 
More generally, if $\gcd(n,i,j) = g$, then $C_n(i,j)$
consists of $g$ components, all isomorphic to $C_{n/g}(i/g, j/g)$.

As the drawing of $C_{10}(1,4)$ in Figure~\ref{fig:EZegVoltage} illustrates, two-generator circulant graphs can be drawn symmetrically. 
More precisely, they are always vertex-transitive.
For any $s \in \Z_n$, let $\sigma_s$ be translation by $s$; that is, $\sigma_s(a) = s+a$ for all $a \in \Z_n$. This is the natural left action of the voltage group on the derived graph and it is easily verified to be a graph automorphism of $C_n(i,j)$.  If $a, b$ are vertices of $C_n(i,j)$, then $\sigma_{b-a}(a) = b$. Note that $\sigma_{b-a}$ maps $\{a, a+i\}$ to $\{b, b+i\}$ and $\{a, a+j\}$ to  $\{b, b+j\}$. Thus $C_n(i,j)$ is edge-transitive if and only if there is an automorphism mapping an edge of the form $\{a, a+i\}$ to an edge of the form $\{b, b+j\}$.
Since the reflection given by $\tau_{-1}(a) = -a$ is also an automorphism of $C_n(i,j)$, circulant graphs are edge-transitive if and only if they are arc-transitive.  The classification of all arc-transitive circulant graphs was found independently by Kovacs~\cite{K2004} and Li~\cite{L2005}.
Based on this work, Poto\u{c}nik and Wilson recently noted in~\cite{PW2020} that $4$-regular, two-generator circulant graphs are edge-transitive if and only if they are either isomorphic to $C_n(1,j)$ for $j^2 \equiv \pm 1$ or isomorphic to $C_{2m}(1, m-1)$ for $m \ge 3$. Thus, for example, $C_{10}(1, 4)$ in Figure~\ref{fig:EZegVoltage} is edge-transitive. The only 3-regular edge-transitive circulant graphs are $C_4(1,2) \cong K_4$ and $C_6(1,3) \cong K_{3,3}$; see~\cite{GKLV2017}.

Another way to characterize the symmetry of a graph $G$ is to compute parameters that measure how easy it is to `break' any nontrivial automorphisms of $G$. As one example of this, a determining set of graph $G$ is a vertex subset $W$ such that the only graph automorphism that fixes each vertex in $W$ is the identity. The size of a minimum determining set is the determining number of $G$, denoted by $\det(G)$. (Some authors refer to this as the fixing number of the graph.) Another example is to assign $d$ colors to the vertices in such a way that the only automorphism that preserves the color classes (setwise) is the identity. Such a coloring is called a $d$-distinguishing coloring; the minimum number of colors required for a distinguishing coloring is called the distinguishing number of the graph, denoted by $\dist(G)$. For more background on determining and distinguishing number and the relationships between them, see~\cite{AB2007}. It has been shown that many infinite families of graphs have distinguishing number $2$; for such graphs, a further refinement is to determine the minimum size of a color class in a $2$-distinguishing coloring. This parameter, introduced in~\cite{B2008}, is called the cost of $2$-distinguishing $G$ and is denoted $\rho(G)$.

If a graph $G$ is disconnected with components $C_1, C_2, \dots C_k$, then it is possible to calculate its symmetry parameters from those of its components. In particular, if all components have positive determining number, then 
$\det(G) = \det(C_1) + \dots + \det(C_k)$.
However, the situation for distinguishing number is more complicated. If multiple components are isomorphic, then there are nontrivial automorphisms that permute components. We therefore need to know the number of nonisomorphic distinguishing colorings for each such component. In this paper, we focus on finding the symmetry parameters only for two-generator circulant graphs that are connected.

Partial results on the determining and distinguishing number of circulant graphs have been obtained. Recently, Brooks et al.~\cite{Brooks2021} studied the determining number of powers of cycles. This motivated their study of general circulant graphs of the form $C_n(A)$, where $A \subseteq \Z_n$ and vertices $u$ and $v$ are adjacent if and only if $\pm(u-v) \in A$. They identify the determining number of circulant graphs with two generators $\{i,j\}$ with $i+j = \frac{n}{2}$, with $i = 1$ and $4 \le j \le \frac{n}{2}$ and, for even $n$, for $i = 2$ and $j>1$ odd.  Brooks et al. conjecture that if $C_n(i,j)$ is connected, then $\det(C_n(i,j)) = 2$ if and only if $C_n(i,j)$ is twin-free. We prove that this is true except for $C_{10}(1,3)$.

Gravier, Meslem and Souad~\cite{GMS2014} investigated the distinguishing number of circulant graphs $C_n(A)$ where $n= mp\ge 3$, for some $m\ge 1$ and $p \ge 2$, and 
$A = \{kp+1 \mid 0 \le k \le m-1\}$.
Restricted to two-generator circulant graphs, their results are
$\dist(C_{2p}(1, p-1)) = 3$, if $p \ge 2$ and $p \neq 4$, and 
$5$, if $p=4$.

The presence of \emph{twin} vertices, which are vertices having the same neighborhood, understandably affects symmetry parameters. For example, Gonzales and Puertas~\cite{GP2019} looked at quotient graphs with respect to the twin relation to find upper and lower bounds on the determining number of an arbitrary graph.
Brooks et al.~\cite{Brooks2021} prove that if every vertex in $C_n(A)$ is in a set of $k$ mutual twins, then $\det(C_n(A)) = n - (n/k)$.

In this paper, we give complete results on the symmetry parameters of connected, two-generator circulant graphs. We begin by characterizing the automorphism group of such graphs. There are results on the automorphisms of special cases of $C_n(A)$, such as when it is arc-transitive or when  $n$ is prime, a prime power, or square-free, and/or the elements of $A$ are divisors of $n$, or when the circulant graph has a rational spectrum; see~\cite{L2005}, \cite{BI2011}, \cite{KK2012} and~\cite{M2007}.  We find the automorphism group of all connected, two-generator circulant graphs, with no restrictions on arc-transitivity or the prime factorization of $n$. As in~\cite{BI2011}, our proofs make extensive use of possible sets of common neighbors. Let $H$ be the set of units in $\mathbb Z_n$ that preserve $\{\pm i, \pm j\}$ under multiplication, commonly denoted by $\aut(\mathbb Z_n, \{\pm i, \pm j\})$. We show that if $C_n(i,j)$ is connected, twin-free, and not $C_{10}(1,3)$, 
then $\aut(C_n(i,j)) = \Z_n \rtimes H$. If $C_n(i,j)$ has twins, then every automorphism can be expressed as the composition of an element of $\Z_n \rtimes H$ and an automorphism that permutes sets of mutually twin vertices.

\begin{figure}[h]
\includegraphics[width=0.65\textwidth, center]{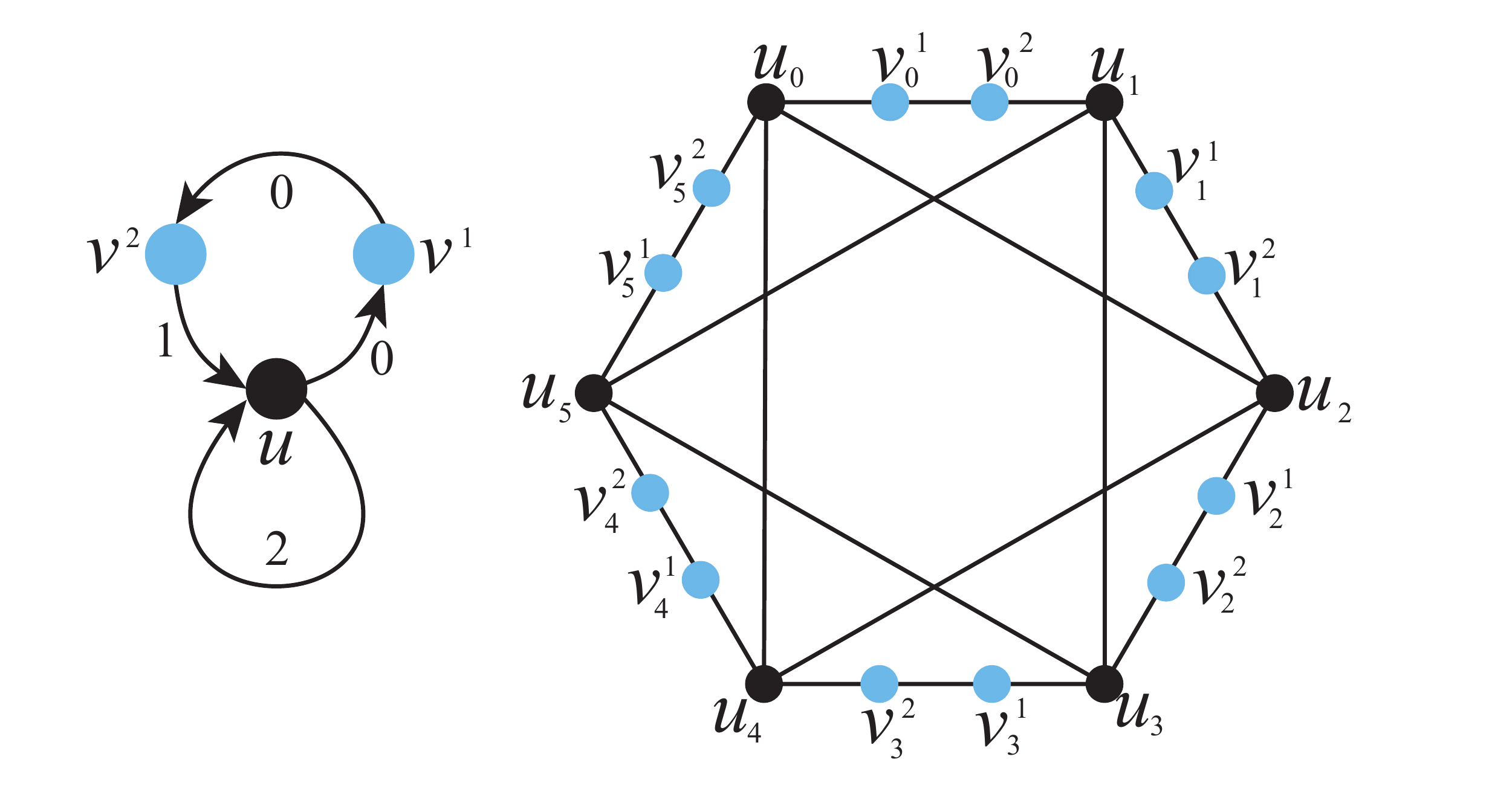}
\caption{Voltage graph subdivided $(B_2, \Z_{6})$ and derived graph $C_6(1_\div 2, 2)$.}\label{fig:SubdividedVoltage}
\end{figure}

 We also consider the derived graphs associated to voltage graphs obtained by subdividing one loop in $B_2$ with $\ell$ vertices of degree $2$. 
 Equivalently, the voltage graph is a directed cycle of length $\ell$ with a loop at one vertex. Using the fact that without loss of generality, we can assign a voltage of $0$ to the arcs in a spanning tree of the base directed graph, we assign the subdivided loop's original voltage to the arc of the cycle directed to vertex $u$~\cite{GT1987}. The associated derived graph is a subdivision of $C_n(i,j)$, which we denote by $C_n(i_\div \ell, j)$ if the arc of voltage $i$ is the one that has been subdivided; $C_n(i, j_\div \ell)$ is defined analogously. See Figure~\ref{fig:SubdividedVoltage}, in which the arc of voltage $i$ in $B_2$ has been subdivided with $\ell=2$ vertices of degree $2$, producing derived graph $C_6(1_\div 2, 2)$. 

We find the automorphism groups of connected, two-generator circulant graphs that have been subdivided in this way, and use them to determine their symmetry parameters. Our results, summarized in Table~\ref{tab:summary}, confirm and extend those in~\cite{Brooks2021} and~\cite{GMS2014}. Note that subdividing an arc of voltage less than $n/2$ sometimes reduces the determining number and always reduces the cost of $2$-distinguishing, indicating that the overall symmetry has been reduced. On the other hand, when an arc of voltage $n/2$ is subdivided, the derived graph changes from trivalent to tetravalent. As a result, the derived graph has additional symmetries and both the determining number and the cost of $2$-distinguishing increase.

\begin{table}[h]
 \centering
 \begin{tabular}{|c|c|c|c|c|}
 \hline
 & $\det$ & $\dist$ & $\rho$ & condition(s)
 \\ \hline
 $C_n(i,j)$ & $n-1$ & $n$ & n/a & $n \in \{4, 5\}$ \\ \cline{2-5}
 &4 & 4 & n/a & $(n,i,j) = (6, 1, 3)$\\
 \cline{2-5}
 & 6 & 5 & n/a & $(n, i, j) = (8,1,3)$\\
 \cline{2-5}
 & 4 & 3 & n/a & $(n,i,j) = (10, 1, 3)$
 \\
 \cline{2-5}
 & $n/2$ & 3 & n/a & $i+j = n/2,$ but $n \neq 8$ \\
 \cline{2-5}
 & 2& 2 & 3 & otherwise
 \\ \hline
 $C_n(i_\div \ell, j)$, or & 1 & 2 & 1 & $\ell\ge 2$ \text{ and } $H = \{\pm 1\}$ \\ \cline{2-5}
 $C_n(i, j_\div \ell),\, j<n/2$ & 2& 2 & 2 & otherwise
 \\ \hline
 $C_n(i, j_\div \ell), \, j=n/2$ & 4 & 3 & n/a & $\ell=1,\, j=2$ \\ \cline{2-5}
 & $j+1$ & 2 & $j+3$ & $\ell=1, \, j\ge 3$ \\ \cline{2-5}
 & $j$ & 2 & $j+1$ & $\ell=2, \, j \in \{2, 3, 4, 5\}$ or $\ell=j=3$ \\ \cline{2-5}
 & $j$ & 2 & $j$ & otherwise
 \\ \hline
 \end{tabular}
 \caption{Summary of Symmetry Parameters}
 \label{tab:summary}
\end{table}

Throughout the entire paper, we assume $
0 < i < j \le n/2$ and $\gcd(n,i,j) = 1$. Thus $C_n(i,j)$ is a connected, two-generator circulant graph.
In Section~\ref{sec:whenTwins}, we characterize which $C_n(i,j)$ have twin vertices. In Section~\ref{sec:TwinCase}, we compute the symmetry parameters of such graphs.
Section~\ref{sec:CommonNeighbors} exhibits the possible sets of common neighbors in twin-free $C_n(i,j)$. In Section~\ref{sec:Automorphisms}, we characterize the automorphisms of $C_n(i,j)$, both for those that are twin-free and, because we use this information in the subdivided case, for those with twins.
Section~\ref{sec:TwinFreeCase} gives the symmetry parameters of twin-free $C_n(i,j)$. Finally, in Section~\ref{sec:SubdividedCase}, we find the automorphism group and symmetry parameters of subdivided $C_n(i,j)$. Subsection~\ref{subsec:jlessn/2} deals with the case in which the subdivided loop in the voltage graph has voltage less than $n/2$, and Subsection~\ref{subsec:jisn/2} considers the case where this voltage equals $n/2$. We close with some ideas for future research in Section~\ref{sec:Open}.

\section{Twins in Two-Generator Circulant Graphs}\label{sec:whenTwins}

The open neighborhood of a vertex $v$ in a graph $G$, $N(v)$, is defined to be $\{u \in V(G) \mid uv \in E(G)\}$ and the closed neighborhood of $v$ is $N[v]= \{v\} \cup N(v)$. Two vertices $v$ and $w$ are nonadjacent twins if $N(v) = N(w)$, and they are adjacent twins if $N[v] = N[w]$. Note that if a vertex has an adjacent twin, it cannot also have a nonadjacent twin. We say a graph has twins if it has adjacent or nonadjacent twins. Finally, we define vertices $u$ and $v$ in a graph $G$ to be \emph{co-twins} if $N[v] = \overline{N[v]}$. 

If $u$ and $v$ are twins, the vertex map that interchanges $u$ and $v$ while leaving all other vertices fixed is a graph automorphism. In circulant graphs, exchanging a pair of co-twins also leads to additional automorphisms. The presence of twin and co-twin vertices therefore has a substantial effect on the symmetry parameters of the graph. 

It is easy to verify that the only connected, two-generator circulant graph with co-twins is $C_{10}(1,3)$; in this case, $N[a] = \overline{N[a+5]}$ for all $a \in \mathbb Z_{10}$. More generally, if $n = 4k+2$ for some $2 \le k \in \mathbb Z$, then  any two vertices of the form $a, a+2k+1$ are co-twins in $C_{2k}(1, 3, \dots, 2k-1)$.

In \cite{Brooks2021}, Brooks et al. show that if $i+j=n/2$, then $N(a) = N(a + n/2)$. Lemma~\ref{lem:whenTwins} strengthens this result by fully characterizing all twin vertices in connected, two-generator circulant graphs.

\begin{lemma}\label{lem:whenTwins}
\begin{enumerate}[(1)]
\item If $n\in \{4, 5\}$, then $C_n(i,j) = C_n(1,2) = K_n$ and so any two distinct vertices are adjacent twins.
\item If $n \ge 6$ and $j <n/2$, then $C_n(i,j)$ is twin-free if and only if $i+j \neq n/2$. If $i+j = n/2$, then for distinct $a,b$, $N(a)=N(b)$ if and only if either $b = a \pm n/2$ or $(n, i, j) = (8, 1, 3)$ and $b = a\pm2$.
\item If $n \ge 6$ and $j = n/2$, then $C_n(i,j)$ is twin-free except if $(n, i, j) = (6,1,3)$, in which case for distinct $a,b$, $N(a) = N(b)$ if and only if $b = a\pm2$.
\end{enumerate}
\end{lemma}

\begin{proof}
Recalling, $0 < i < j \le n/2$, statement (1) is easy to verify. So assume $n \ge 6$.

Since $C_n(i,j)$ is vertex-transitive, $C_n(i,j)$ has twin vertices if and only if $0$ has a twin. Suppose vertex $0 \not \equiv a \in \Z_n$ is a twin of $0$, so that $N(0) = N(a)$.

\medskip

Because the proof is simpler, we first consider (3). Assume $j = n/2$.
Then $j \equiv -j$ and so
$N(0) = \{i, -i, j\} = \{a+i, a-i, a+j \} = N(a)$.
Since $0 \not \equiv a$, $a+k \not \equiv k$ for $k \in \{i,j,-i\}$. 
There are only two possible values of $a+i$.\\
 \noindent{\bf Case 3a.} If $a+i\equiv-i$, then $a+j \equiv i$ and $a-i \equiv j$. Then $a \equiv -2i \equiv i+j$, so $3i \equiv -j \equiv n/2$. Thus $\gcd(n, i, j) = i$, which by assumption implies $i=1$, and hence $(n, i, j) = (6, 1, 3)$.\\
 \noindent{\bf Case 3b.} If $a+i\equiv j$, then $a-i \equiv i$ and $a+j \equiv -i$. Then $a \equiv 2i \equiv j-i$, so again $3i \equiv j \equiv n/2$, which forces $(n,i, j) = (6, 1, 3)$.\\
Thus, if $N(0)=N(a)$, then $(n,i, j) = (6, 1, 3)$. Moreover, $a\equiv \pm 2 i \equiv \pm 2$. By vertex-transitivity, any two vertices whose difference in $\Z_n$ is $\pm 2$ will be nonadjacent twins. 

 \medskip
 
Returning to (2),
assume $j < n/2$. Then assuming $N(0) = N(a)$ gives $\{i, -i, j ,-j\} = \{a+i, a-i, a+j, a-j\}$.
We again divide into cases based on the possible values of $a+i$.

\medskip 

\noindent {\bf Case 2a.} First, suppose $a+i \equiv -j$, which implies $a \equiv -i -j$ and so $a+j \equiv -i$, There are only two possibilities for $a-i$ and $a-j$.
\begin{itemize}
 \item If $a-i\equiv i$ and $a-j \equiv j$, then $a \equiv 2i \equiv 2j$. Since $0<i <j < n/2$, $0 < 2i < 2j < n$. This is a contradiction. 
 \item If $a-i\equiv j$ and $a-j \equiv i$, then $a \equiv i+j \equiv -i - j$, which implies that $2i+2j\equiv 0$ and so $i+j = n/2$. 
\end{itemize}

\noindent {\bf Case 2b.} Next, suppose $a+i \equiv -i$. 
\begin{itemize}
 \item If $a-i \equiv i$, then because $a+j\not \equiv j$, it must be the case that $a+j\equiv -j$ and $a-j \equiv j$. These imply $a \equiv 2i \equiv -2i \equiv 2j \equiv -2j$. In particular, this implies $2i \equiv 2j$, a contradiction.
 \item Next suppose $a-i \equiv j$. This implies $a-j \equiv i$ and so by process of elimination, $a+j \equiv -j$. These imply $a \equiv -2i \equiv i+j \equiv -2j$. But $-2i \equiv -2j$ implies the contradiction $2i \equiv 2j$.
 \item Last, suppose $a -i \equiv -j$. This implies $a+j \equiv i$ and so $a-j \equiv j$. These imply $a \equiv -2i \equiv i-j \equiv 2j$, which in turn implies $2i+2j \equiv 0$, and so as before, $i + j = n/2$.
\end{itemize}

\noindent {\bf Case 2c.} Finally, suppose $a+i \equiv j$, which implies $a-j \equiv -i$. In this case, there are only two possibilities for $a-i$ and $a+j$. 
\begin{itemize}
 \item If $a-i \equiv -j$ and $a+j \equiv i$, then $a \equiv i-j \equiv -i+j$. This implies $2i \equiv 2j$ which in turn implies $i=j$, a contradiction.
 \item If $a-i\equiv i$ and $a+j \equiv -j$, then $a \equiv 2i \equiv -2j \equiv -i+j$. As argued in Case $2$, this implies $i+j = n/2$.
\end{itemize}

Considering the last option in each of these three cases, we conclude that if $N(0) = N(a)$, then either $a \equiv i+j$ and $i+j = n/2$, or $a \equiv j-i\equiv 2i \equiv -2j$, which indirectly implies $i+j = n/2$. By contrapositive, if $i + j \neq n/2$, then $C_n(i,j)$ is twin-free. If $i+j = n/2$, then it is easy to verify that $N(0)=N(i+j) = \{i, j, -i, -j\}$, so by vertex transitivity, any two vertices of the form $a$ and $a+i+j$ are nonadjacent twins. 
If, in addition, $j-i\equiv 2i \equiv -2j$, then 
the assumption that $\gcd(n, i, j) = 1$ forces $(n, i, j) = (8, 1, 3)$. In this case, $N(2i) = N(2) = \{1, 3, 5, 7\} = N(0)$.
By vertex transitivity, any two vertices of the form $a$ and $a\pm 2$ are nonadjacent twins in this case. 
\end{proof}

\begin{cor}\label{cor:TwinsEdgeTrans}
Every connected, two-generator circulant graph with twins is edge-transitive.
\end{cor}

\begin{proof}

Clearly $C_4(1,2) \cong K_4$, $C_5(1,2) \cong K_5$ and $C_6(1,3) \cong K_{3,3}$ are edge-transitive. By Lemma~\ref{lem:whenTwins}, it suffices to consider $C_n(i,j)$ where $i+j = n/2$. Let $m=n/2$. Since $\gcd(2m,i,j) = 1$ and $2i+2j=2m$, we have $\gcd(i,j) = 1$.

Suppose $j$ is not a unit. Then $\gcd(j,2m)>1$. Furthermore, $\gcd(j,2m)$ divides $2i$, but not $i$. Hence, $j$ is even. If $i$ is also not a unit, we find that $i$ is even, contradicting $\gcd(i,j) = 1$. Thus, without loss of generality, we may assume that $i$ is a unit so $C_n(i,j) \cong C_{2m}(1, i^{-1}j)$. Since $2i+2j \equiv 0$, $2 + 2i^{-1}j \equiv 0$, implying $i^{-1}j \equiv m-1$. 
\end{proof}

The converse of Corollary~\ref{cor:TwinsEdgeTrans} is false; for example, $C_{15}(1,4)$ is edge-transitive and twin-free. Figure~\ref{fig:TwinEgs} shows some examples of circulant graphs with twins.

\begin{figure}
\includegraphics[width= 0.9\textwidth, center]{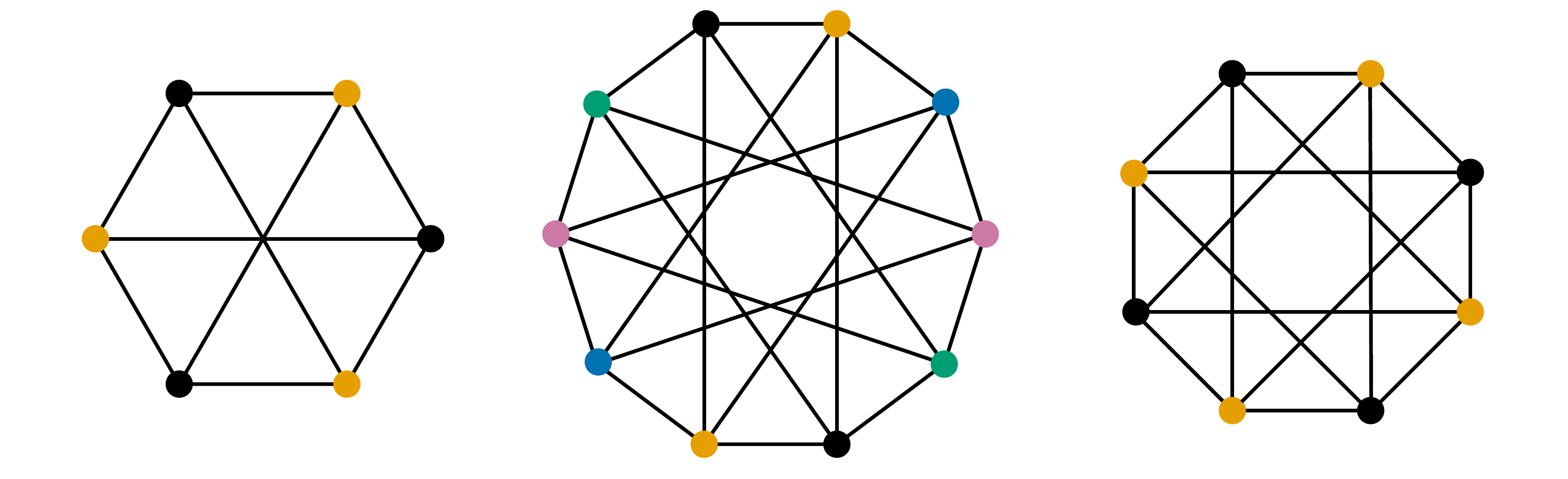}
\caption{Circulant graphs with twins: $C_6(1,3)$, $C_{10}(1,4)$ and $C_8(1,3)$.}
\label{fig:TwinEgs_original_}
\end{figure}

\medskip

\section{Symmetry Parameters for $C_n(i,j)$ with Twins}\label{sec:TwinCase}

As noted at the beginning of Section~\ref{sec:whenTwins}, interchanging a pair of twin vertices and fixing all other vertices is nontrivial graph automorphism. More generally, a vertex map that permutes the elements of a set of mutual twins and fixes all other vertices is a nontrivial graph automorphism.
This implies that a determining set for $G$ must contain at all but one vertex from any set of mutually twin vertices. Additionally, in any distinguishing coloring, the vertices in a set of mutually twins must be assigned distinct colors.

For vertices $x$, $y$ of a graph $G$, define the relation $x \sim y$ if $x$ and $y$ are twin vertices. It is easy to verify that $\sim$ is an equivalence relation on $V(G)$ and so we can create a quotient graph $\widetilde G$ with respect to the relation $\sim$, where the set of vertices is the set equivalence classes $[x] = \{y \in V(G) \mid x \sim y\}$ with $[x]$ adjacent to $[z]$ in $\widetilde G$ if and only if there exist $u \in [x]$ and $v \in [z]$ such that $u$ and $v$ are adjacent in $G$. We refer to $\widetilde G$ as the \emph{twin quotient graph}. It can be verified that $\widetilde G$ is twin-free; see~\cite{BCKLPR2020b}.

Using terminology from~\cite{BCKLPR2020b}, a {\it minimum twin cover} $T$ of a graph $G$ is a subset of vertices that contains all but one vertex from each set of mutual  twins. The following is a corollary of Theorem~19 of~\cite{BCKLPR2020b}: 

\begin{cor}\label{cor:BCKLPR}
Assume every vertex of $G$ has at least one twin. Let $T$ be minimum twin cover of $G$. Then $T$ is a minimum size determining set for $G$. 
\end{cor}

The following result, based on an approach used by Boutin and Cockburn in~\cite{BC2021} to find the symmetry parameters of orthogonality graphs, gives a relationship between the distinguishing numbers of a graph $G$ and $\widetilde G$.

\begin{thm}\label{thm:DistTwins}
Let $G$ be a graph in which every vertex is in a set of exactly $k$ mutual twins and let $\widetilde G$ be the corresponding twin quotient graph. If $\dist(\widetilde G) = \widetilde d$, then $\dist(G) = d$, where $d$ is the smallest positive integer such that $\binom{d}{k} \ge \widetilde d$.
\end{thm}

\begin{proof}
Assume $\dist(\widetilde G) = \widetilde d$. If $d \in \mathbb N$ satisfies $\binom{d}{k} \ge \widetilde d$, then from a palette of $d$ colors, we can create $\widetilde d$ distinct subsets of $k$ colors, which we can think of a $\widetilde d$ distinct $k$-color packets. By assumption, we can then color each vertex of $\widetilde G$ with a $k$-color packet in a distinguishing fashion. For each $x \in V(G)$, the equivalence class $[x] \in V(\widetilde G)$ has been assigned a $k$-color packet; we randomly assign these $k$ colors to the $k$ vertices in $[x]$ a bijective fashion. 

To show that this defines a distinguishing coloring of $G$, assume $\alpha \in \aut(G)$ preserves the color classes. Since automorphisms preserve neighborhoods and vertices of $G$ are identified in $\widetilde G$ precisely when they have the same neighbors (except possibly each other), $\alpha$ induces an automorphism $\widetilde \alpha$ of $\widetilde G$ given by $\widetilde \alpha([x]) = [\alpha(x)]$. Since $\alpha$ preserves color classes in $G$, the same $k$ colors appear in the $k$ vertices of $[x]$ and the $k$ vertices of $\widetilde \alpha([x])$, which means that $\widetilde \alpha$ preserves the $k$-color packets assigned to the vertices of $\widetilde G$. By assumption, 
$\widetilde \alpha$ is the identity automorphism on $\widetilde G$. Hence, $\alpha$ is a coloring-preserving bijection on each equivalence class of vertices. However, the vertices of each equivalence class have distinct colors, and so $\alpha$ must be the identity automorphism on $G$. This shows $\dist(G) \le d$.

Next, let $d'<d$. Then $\binom{d'}{k} < \widetilde d$. Assume there exists a $d'$-distinguishing coloring of $G$. The $k$ vertices in each equivalence class of mutual twins must be assigned $k$ different colors. Thus we can assign a $k$-color packet to each vertex of $\widetilde G$. Since the number of different $k$-color packets possible is less than $\widetilde d \dist(\widetilde G)$, this coloring cannot be distinguishing, and so there exists a nontrivial automorphism $\widetilde \gamma$ of $\widetilde G$ that preserves the coloring with $k$-color packets. 

We can define a corresponding automorphism $\gamma$ of $G$ as follows. For each $x \in V(G)$, the $k$-color packet assigned to $[x]$ is the same as the $k$-color packet assigned to $\widetilde \gamma([x])$, by assumption. For each $y \in [x]$, we let $\gamma(y)$ be the vertex in $\widetilde \gamma([x])$ that has the same color as $y$ in the $d'$-distinguishing coloring of $G$. It is easy to verify that $\gamma$ is a nontrivial, color-preserving automorphism of $G$, a contradiction. Hence $\dist(G) \ge d$, and we are done.
\end{proof}

\begin{thm}\label{thm:DetDistTwins}
Assume $C_n(i,j)$ has twins.
\begin{enumerate}[(1)]
\item If $n\in \{4, 5\}$, then $\det(C_n(i,j)) = n-1$ and $\dist(C_n(i,j))= n$.
\item If $(n, i, j) = (6, 1,3)$, then $\det(C_n(i,j)) = \dist(C_n(i,j))= 4$.
\item If $n \ge 6$ and $i+j=n/2$ then $\det(C_8(1,3)) = 6$ and $\dist(C_8(1,3))=5$ and for all other values, 
$\det(C_n(i,j)) = n/2$ and $\dist(C_n(i,j)) = 3$.
\end{enumerate}
\end{thm}

\begin{proof}
The three statements in this theorem align with the three cases in Lemma~\ref{lem:whenTwins}. Statement (1) handles the cases in which $C_n(i,j) = K_n$. 

For statement (2), we know from Lemma~\ref{lem:whenTwins} that vertices $0, 2, 4$ are mutually nonadjacent twins, as are $1, 3, 5$. A minimum twin cover is $T= \{0, 1, 2, 3\}$; since every vertex has a nonadjacent twin, by Corollary~\ref{cor:BCKLPR} this is also a minimum determining set. The twin quotient graph in this case is $K_2$ which has distinguishing number $\widetilde d = 2$ and each equivalence class contains $k=3$ vertices. In this case, we need $d=4$ to get $\binom{d}{3} \ge 2$, so $\dist(C_n(i,j)) = 4$ by Theorem~\ref{thm:DistTwins}. 

For statement (3), we first consider the case $(n,i,j) = (8, 1, 3)$. By Lemma~\ref{lem:whenTwins}, vertices $0, 2, 4, 6$ are mutually nonadjacent twins, as are vertices $1, 3, 5, 7$. A minimum twin cover is $T = \{0, 1, 2, 3, 4, 5\}$ and so by Corollary~\ref{cor:BCKLPR} this is also a minimum determining set. The twin quotient graph is again $K_2$ which has determining number $\widetilde d = 2$; this time each equivalence class contains $k=4$ vertices. We need $d=5$ to get $\binom{d}{4} \ge 2$, so $\dist(C_n(i,j)) = 5$ by Theorem~\ref{thm:DistTwins}. 

Finally, we consider the case $n \ge 6$, $i+j=n/2$ but $n \neq 8$. From Lemma~\ref{lem:whenTwins}, $a$ and $a+i+j=a+n/2$ are nonajdacent twins for all $a \in \Z_n$. A minimum twin cover and minimum determining set is $T = \{0,1,2,\ldots,n/2-1\}$. The twin quotient graph has order $n/2$ and is connected. Since $N(a) =\{a+i,a_i,a+j,a-j\}$, in the twin quotient graph we have $N(\{a, a+i+j\}) = \lbrace \{a+i, a-i\}, \{a+j, a-j\}\rbrace$. Hence the twin quotient graph is 2-regular. It follows that it is $C_{n/2}$, which has distinguishing number $\tilde d = 3$ if $n/2 \in \{3, 4, 5\}$ and distinguishing number $\tilde d = 2$ if $n/2 \ge 6$. In either case, the smallest $d$ such that $\binom{d}{2} \ge \tilde d$ is $d=3$. 
\end{proof}

Since there is only one two-generator graph with co-twins, namely $C_{10}(1,3)$, we  use direct computation to find the symmetry parameters, $\det(C_{10}(1, 3)) = 4$ and $\dist(C_{10}(1, 3)) = 3$.

\section{Common Neighbors in Twin-free $C_n(i,j)$}\label{sec:CommonNeighbors}

To find the determining and distinguishing numbers of $C_n(i,j)$ in the twin-free case, we first find its automorphism group. A key tool in our investigation involves possible sets of common neighbors. Since automorphisms respect adjacency and nonadjacency, if $\alpha \in \aut(G)$ and $u, v$ are vertices in a graph $G$ with 
$
N(u) \cap N(v) = \{w_1, \dots w_\ell\},
$
then 
$
N(\alpha(u)) \cap N(\alpha(v)) = \{\alpha(w_1), \dots \alpha(w_\ell)\}$.
To use this fact, we must determine the possible sets of common neighbors.

Two vertices have common neighbors if and only if there is a path of length 2 between them. In $C_n(i,j)$, we find $a$ and $b$ have common neighbors if and only if 
$b \in \{a+2i, a-2i, a+2j, a-2j, a+i+j, a+i-j, a-i+j, a-i-j\}$.
Exchanging the roles of $a$ and $b$ as needed, we need only consider the cases $b \in \{a+2i, a+2j, a+i+j, a+i-j\}$.

\begin{lemma}\label{lem:CommonNeighborjisn/2}
Assume $C_n(i,j)$ is twin-free and $j = n/2$.
Then distinct vertices $a, b$ have common neighbors if and only if $b \in \{a+2i, a+i+j\}$. Furthermore,
\begin{enumerate}[(1)]
\item $N(a) \cap N(a+2i) = \{a+i\}$,
\item $N(a) \cap N(a+i+j) = \{a+i, a+j\}$.
\end{enumerate}
\end{lemma}

\begin{proof}
By Lemma~\ref{lem:whenTwins}, the assumption that $C_n(i,j)$ is twin-free implies $n \ge 6$ 
and $(n,i,j) \neq (6,1,3)$.
The assumption that $j=n/2$ implies 
$j \equiv -j $ and $2j \equiv 0$, so 
without loss of generality 
$a$ and $b$ have common neighbors only when $b \in \{a+2i, a+i+j\}$.

By vertex transitivity, it suffices to consider the case $a=0$. Note that 
$
N(0)=\{i, -i, j\}$, $N(2i) = \{i, 3i, 2i+j\}$ and $N(i+j) = \{i, j, 2i+j\}$.
For (1), clearly $i \in N(0) \cap N(2i)$.
We show the other candidates lead to a contradiction.
\begin{itemize}
 \item If $-i \equiv 3i$, then $4i\equiv 0$. Since we assumed $0<i<j=n/2$, this requires $2i=j$. But then $i=\gcd(n,i,j)=1$ requires $i=1$ and $n=4$, a contradiction. 
 \item If $-i \equiv 2i+j$ we get $3i \equiv -j \equiv j$, combining two cases. If $3i \equiv j$, then $i = \gcd(n,i,j) = 1$ implies $(n,i,j) = (6,1,3)$, a contradiction.
 \item If $j \equiv 2i+j $, then $2i \equiv 0$, contradicting $0<i < n/2$.
\end{itemize}

For (2), clearly $\{i, j \} \subseteq N(0) \cap N(i+j)$. Since all vertices have degree 3 and $C_n(i,j)$ is twin-free, 
they 
cannot have a third common neighbor. 
\end{proof}

\begin{lemma}\label{lem:CommonNeighbors}
Assume $C_n(i,j)$ is twin-free and $j < n/2$. 
Then distinct vertices $a, b$ have common neighbors if and only if $b \in \{a+2i, a+2j, a+i+j, a+i-j\}$. Furthermore, the common neighbors of each such pair of vertices $(a,b)$ is given by the following table, where we take the union of the applicable rows.

\begin{center}
\begin{tabular}{|c|c|c|c|c|}
\hline
 & $(a,a+2i)$ & $(a,a+2j)$ & $(a,a+i+j)$ & $(a,a+i-j)$\\ \hline 
 always & $\{a+i\}$ & $\{a+j\}$ & $\{a+i, a+j\}$ & $\{a+i, a-j\}$ \\ \hline
 $4i\equiv0$ & $\{a-i\}$ & --- & --- & ---\\
 $4j\equiv0$ & --- & $\{a-j\}$ & --- & ---\\ \hline
 $3i\equiv-j$ & $\{a-i, a-j\}$ & --- & $\{a-i\}$ & ---\\
 $3i\equiv j$ & $\{a-i, a+j\}$ & --- & --- & $\{a-i\}$\\ \hline 
 $3j\equiv -i$ & --- & $\{a-i, a-j\}$ & $\{a-j\}$ & ---\\
 $3j\equiv i$ & --- & $\{a+i, a-j\}$ & --- & $\{a+j\}$\\ \hline
\end{tabular}\label{tab:CommonNeighbors}
\end{center}

\end{lemma}

\begin{proof} By Lemma~\ref{lem:whenTwins}, the assumption that $C_n(i,j)$ is twin-free implies $n \ge 6$ and $i+j \neq n/2$. 
Again, by vertex transitivity, it suffices to consider the case $a=0$. As in the proof of the previous lemma, we will use contradictions to rule out possible common neighbors. The assumption $0 < i < j < n/2$ implies 
$2i \not\equiv 0$, $2j \not\equiv 0 $ and $2i \not\equiv 2j$.
The assumption that $i + j \neq n/2$ implies 
$2i \not\equiv -2j$. 
Note that $i$ and $j$ are interchangeable in these contradictions.

\medskip

For $(a,a+2i)$, we start with
$N(0)=\{i, -i, j, -j\}$ and $N(2i) = \{3i, i, 2i+j, 2i-j\}$.
Clearly $i \in N(0) \cap N(2i)$. There are three ways $-i$ could also be a common neighbor.
\begin{itemize}
 \item If $-i \equiv 3i $, then $4i \equiv 0 $. 
 In this case, there are no other common neighbors because all other options require $2i \equiv 0$ or $2i \equiv \pm 2j$.
 \item If $-i \equiv 2i+j$, then $3i \equiv -j$, so there are three common neighbors: $i$, $-i$ and $-j$. Because $C_n(i,j)$ is twin-free, $j$ cannot be a fourth common neighbor. 
 \item If $-i \equiv 2i-j$, then $3i \equiv j $, so again there are three common neighbors: $i$, $-i$ and $j$. 
\end{itemize}

Next assume that $-i \notin N(0) \cap (N(2i)$. We will show that in this case, neither $j$ nor $-j$ can be a common neighbor. If $j \equiv 3i$, then $2i-j \equiv -i$, meaning that $-i$ is a common neighbor, a contradiction. If $j \equiv 2i+j$ then $2i \equiv 0$, and if $j \equiv 2i-j$, then $2i \equiv 2j$; both contradictions. Similarly, assuming $-j \equiv 3i$ implies that $-i \equiv 2i-j$ is a common neighbor; assuming $-j \equiv 2i+j$ and $-j \equiv 2i-j$ lead to the contradictions $2i+2j \equiv 0$ and $2i\equiv 0$ respectively.

\medskip

For $(a,a+2j)$, the argument is analogous to the one for $(a,a+2i)$, with $j$ replacing $i$. 

\medskip

For $(a,a+i+j)$, we start with 
$N(0)=\{i, -i, j, -j\}$ and $N(i+j) = \{2i+j, j, i+2j, i\}$.
Clearly $\{i, j\} \subseteq N(0) \cap N(i+j)$.

We consider the options for an additional common neighbor. If $-i \equiv 2i+j$ or $-j \equiv i+2j$, then $i+3j \equiv 0$ or $3i+j \equiv 0$, respectively. 
In each, the assumption that $C_n(i,j)$ is twin-free means there are no additional common neighbors. 
In the remaining cases, $-i \equiv i+2j$ and $-j \equiv 2i+j$, we have $-2i \equiv 2j$, a contradiction. 

\medskip

For $(a,a+i-j)$, the argument is analogous to that for $(a,a+i+j)$. 

\end{proof}

It is possible for two special conditions to hold, in which case both affect the set of common neighbors. For example, in $C_{12}(3, 5)$, both conditions $4i\equiv 0$ and $3j\equiv i$ hold. In this case, $N(0)\cap N(2i) = \{i, -i\}$, $N(0) \cap N(i+j) = \{i, j\}$ and 
$N(0)\cap N(2j) =N(0)\cap N(i-j)= \{i, j, -j\}$.

\begin{figure}[h]
\includegraphics[width=0.4\textwidth, center]{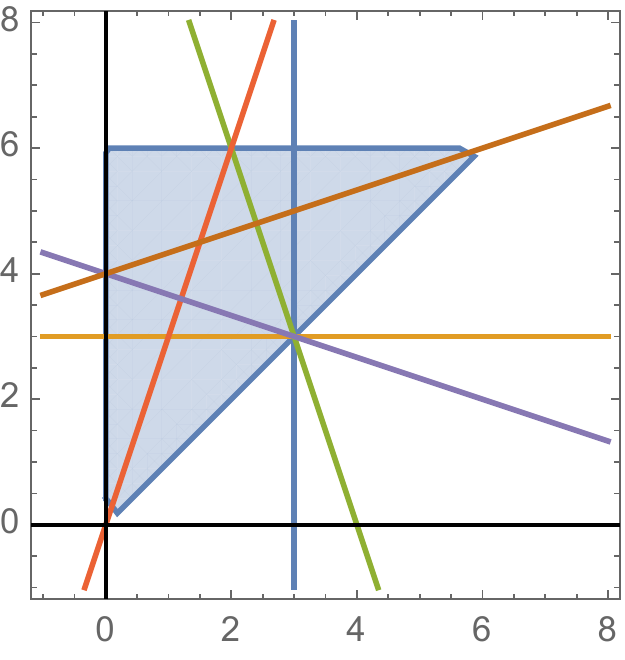}
\caption{Lines corresponding to special conditions when $n=12$.}\label{fig:lines}
\end{figure}

Our assumption that $0<i<j< n/2$ implies that each condition corresponds to exactly one linear equation in $\mathbb{Z}$. We can investigate when multiple special conditions hold simultaneously by graphing these equations within this region. Figure~\ref{fig:lines} illustrates these lines in the $ij$-plane when $n=12$. Observe that there are five intersection points, each of exactly two lines. Given that we are interested only in intersection points corresponding to integral values of $n$, $i$ and $j$, that also satisfy $\gcd(n, i, j) = 1$ and $i+j \neq n/2$, we find only three cases where multiple special conditions are satisfied: 
$C_{12}(3,5)$ satisfies both $4i \equiv 0$ and $3j\equiv i$, $C_{10}(1,3)$ satisfies both $3i\equiv j$ and $3j\equiv -i$ and $C_{12}(1,3)$ satisfies both $4j \equiv 0$ and $3i\equiv j$.

\section{Automorphisms of Circulant Graphs}\label{sec:Automorphisms}

In the introduction, we noted that for any $s\in \Z_n$, translation by $s$, given by $\sigma_s(a) = s+a$, is a graph automorphism of $C_n(i, j)$. Thus $\aut(C_n(i,j))$ has a subgroup isomorphic to $\Z_n$. Additionally, let $H$ be the set of automorphisms $\mathbb Z_n$ that preserve $\{\pm i, \pm j\}$ ; that is,  $H=\aut(\mathbb Z_n, \{\pm i, \pm j\})$.
It is easy to verify that $H$ is a subgroup of $\aut \mathbb(\mathbb Z_n) = U(n)$. 
For $t\in H$ and $a \in \Z_n$, let $\tau_t(a) = ta$. The following is a special case of Godsil's Lemma~2.1 in~\cite{G1981}. For the sake of completeness, we include a proof.

\begin{prop}
 For all $t \in H$, we have $\tau_t\in \aut(C_n(i,j))$.
\end{prop}

\begin{proof}
 By definition, vertices $a$ and $b$ in $C_n(i,j)$ are adjacent if and only if $a-b \equiv \pm i$ or $\pm j$. Since $\tau_t(a) - \tau_t(b) = ta - tb = t(a-b)$, $a$ and $b$ are adjacent if and only if $\tau_t(a)$ and $\tau_t(b)$ are adjacent. 
\end{proof}

It will always be the case that $\{\pm 1\} \subseteq H$. In fact, in many cases, these are the only two elements of $H$. Lemma~\ref{lem:HedgeTrans}  gives $H$ in the edge-transitive cases.

\begin{lemma}\label{lem:HedgeTrans}
For $C_4(1,2)$, $C_6(1,3)$ and $C_{2m}(1, m-1)$ where $m \ge 3$ is odd, $H=\{\pm 1\}$. 
For $C_{2m}(1, m-1)$ where $m \ge 3$ is even and $j = m-1$ and $C_n(1,j)$ where $j^2 \equiv \pm 1$, $H = \{\pm 1, \pm j\}$.
\end{lemma}

\begin{proof}
The results for $C_4(1,2)$ and $C_6(1,3)$ follow from $U(4) = U(6) = \{ \pm 1 \}$. Since $t \in H$ preserves $\{\pm 1, \pm j\}$, $t=t\cdot 1 \in \{\pm 1, \pm j\}$, so
$H \subseteq \{\pm 1, \pm j\}$.
If $m$ is odd, then $j = m-1$ is even and so $\pm j \notin U(2m) = U(n)$. If $m$ is even, then $2m$ divides $m^2$ and so $j^2 =(m-1)^2 = m^2 - 2m+1 \equiv 1$. In particular, $j \in U(n)$.
If $j^2 \equiv \pm 1$, then clearly $j$ and $-j$ preserve $\{\pm 1, \pm j\}$.
\end{proof}

If $C_n(i,j)$ is not edge-transitive, then as noted in the introduction, no automorphism of $C_n(i,j)$ takes an edge of the form $\{a, a+i\}$ to an edge of the form $\{b, b+j\}$. Hence in this case, $H = \aut(\mathbb Z_n, \{\pm i\}, \{\pm j\})$.

\begin{lemma}\label{lem:His1-1}
Let $0<i<j \le n/2$ and $\gcd(n,i,j) = 1$. If $t\in U(n)$ satisfies $ti \equiv i$ and $tj \equiv j$, then $t\equiv 1$. If $t$ satisfies $ti \equiv -i$ and $tj \equiv -j$, then $t\equiv -1$.
\end{lemma} 

\begin{proof}

Since $\gcd(n,i,j) = 1$, there exist $x, y, z \in \Z$ such that $xi+yj+zn = 1$, which means $xi+yj \equiv 1$. Hence 
\[
t \equiv t(xi+yj) \equiv x(ti) + y(tj) \equiv \begin{cases}
xi+yj\equiv 1, \,&\text{ if } ti \equiv i \text{ and } tj \equiv j,\\
-(xi + yj)\equiv -1, &\text{ if } ti \equiv -i \text{ and } tj \equiv -j.
\end{cases}
\]
\end{proof}

\begin{cor}\label{cor:Hnot1-1}
Assume $C_n(i,j)$ is not edge-transitive. If there exists $1 \neq t \in U(n)$ such that $ti \equiv i$ and $tj \equiv -j$, then $H = \{\pm 1, \pm t\}$. Otherwise $H = \{\pm 1\}$.
\end{cor}

\begin{proof}
If there exists such a $t$, then  $\pm t\in H= \aut(\mathbb Z_n, \{\pm i\}, \{\pm j\})$. 
Suppose $t^* \in H \setminus \{ \pm 1, \pm t\}$. 
By Lemma~\ref{lem:His1-1}, we can assume without loss of generality that $t^*i \equiv i$ and $t^*j \equiv -j$. Then $tt^*$ and $t^2$ fix both $i$ and $j$ and so by Lemma~\ref{lem:His1-1}, $tt^* \equiv t^2 \equiv 1.$ Since $t$ is a unit, this implies $t\equiv t^*$.
\end{proof}

\begin{figure}
\includegraphics[width= 0.45\textwidth, center]{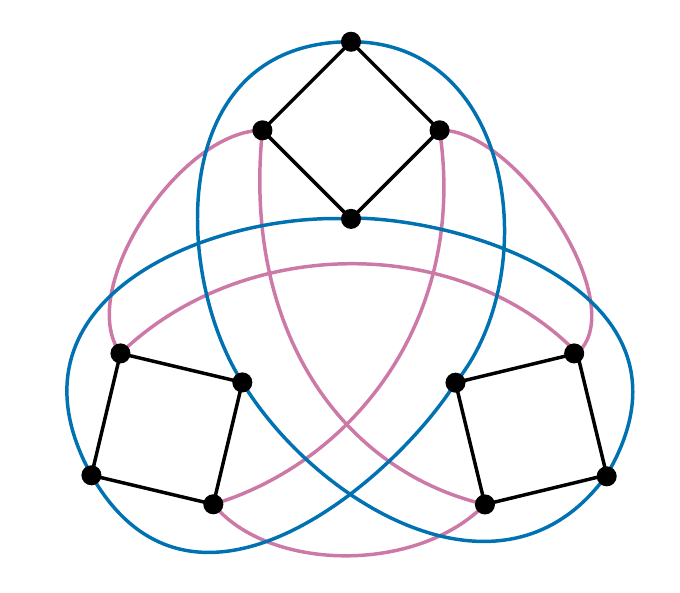}
\caption{$C_{12}(2,3)$.}
\label{fig:ExC12_2_3}
\end{figure}

\begin{Ex}
 For a non-edge-transitive example where $H \neq \{\pm 1\}$, let $(n, i, j) = (12, 2, 3)$. 
 Then $U(12)=\{1, 5, 7, 11\}$. 
 Note that $t = 7$ satisfies $ti \equiv +i$ and $tj \equiv -j$.
 The (twin-free) circulant graph $C_{12}(2,3)$ is shown in Figure~\ref{fig:ExC12_2_3}. 
\end{Ex} 

The possibilities for $H$ are summarized in Table~\ref{tab:PossibleH}.

\begin{table}[h]
    \centering
    \begin{tabular}{|c|c|c|}
    \hline
    & $H $ & Conditions \\
    \hline
     & $\{ \pm 1, \pm j\}$ & $C_{2m}(1, m-1)$,\, $m \ge 3 $ even ($j = m-1$)\\
    Edge-transitive & & $C_n(1, j), \, j^2 \equiv \pm 1$\\
    \cline{2-3}
    $H = \aut(\mathbb Z_n, \{\pm i, \pm j\})$ & $\{\pm 1\}$ & $C_{2m}(1, m-1)$,\, $m \ge 3$ odd\\
    & & $C_4(1,2), \, C_6(1,3)$\\
    \hline
    Not edge-transitive & $\{\pm 1, \pm t\}$ & $  1 \neq t \in U(n)$ satisfies $ti \equiv i$, $tj \equiv -j$ \\
    \cline{2-3}
    $H = \aut(\mathbb Z_n, \{\pm i\},\{ \pm j\})$& $\{ \pm 1\}$ & otherwise\\
    \hline    
    \end{tabular}
    \caption{Possibilities for $H$}
    \label{tab:PossibleH}
\end{table}

For any $s \in \Z_n$ and $t \in H$, we can compose the automorphisms $\sigma_s$ and $\tau_t$; for any $a \in \Z_n$,
$(\sigma_s \circ \tau_t) \cdot (a) = s+ta$.
These automorphisms do not commute, but
$(\tau_t)^{-1} \circ \sigma_s \circ \tau_t = \sigma_{st^{-1}}$. In the next subsection, we will show that these are the only automorphisms of connected, twin-free, two-generator circulant graphs, except for $C_{10}(1,3)$.

\subsection{Automorphisms of Twin-Free $C_n(i,j)$}

\begin{prop}\label{prop:autoFix0}
If $C_n(i,j)$ is not $C_{10}(1,3)$, is connected, and twin-free and 
$\alpha \in \aut(C_n(i,j))$ satisfies $\alpha(0) = 0$, then $\alpha$ is an automorphism of the additive group $\Z_n$. 
\end{prop}

\begin{proof}
Assume $\alpha \in \aut(C_n(i,j))$ fixes $0$. Since $\gcd(n, i, j) = 1$, there exist $x, y \in \Z$ such that $xi + yj \equiv 1$. It follows that for any $a \in \Z_n$, there exist $c,d \in \Z$ such that $ci + dj \equiv a$. It suffices to show that for all $0 \le c, d \in \Z$, 
\begin{equation}\label{eqn:induction}
\alpha(ci +dj) \equiv c\alpha(i) + d\alpha(j).
\end{equation}

This proof involves multiple cases, but the underlying strategy is to use induction on $m= c+d$. The base case $m=1$ is trivial; sometimes we need an additional base case. As indicated at the beginning of Section~\ref{sec:TwinCase}, the main tool is to apply $\alpha$ to an equation expressing a set of common neighbors.

First assume $j = n/2$. Since $2j \equiv 0$, we only need to show $\alpha(ci+dj) \equiv c\alpha(i) + dj$ for all $c \ge 0$ and $d \in \{0,1\}$.

Consider a second base case, $m=2$, first with $c=2$ and $d=0$. From Lemma~\ref{lem:CommonNeighborjisn/2}, $N(0) \cap N(2i) = \{i\}$. Applying $\alpha$ gives $N(0) \cap N(\alpha(2i)) = \{\alpha(i)\}$. Since $\alpha(2i)$ has exactly one common neighbor with $0$, we have $\alpha(2i) \in \{2i,-2i\}$ and $\alpha(i) \in \{i,-i\}$ with matching sign. Hence, $\alpha(2i) = 2\alpha(i)$ and equation~\ref{eqn:induction} holds.

We next show equation~\ref{eqn:induction} holds when $c=d=1$. 
Applying $\alpha$ to $N(0) \cap N(-2i) = \{-i\}$ shows $\alpha(-i) \equiv \mp i$ and $\alpha(-2i) \equiv 2\alpha(-i)$. Since $\alpha$ is injective, $\alpha(-i) \equiv -\alpha(i)$ and $\alpha(-2i) \equiv -2 \alpha(i)$.
Next, $\alpha(j)$ must be adjacent to $\alpha(0) = 0$, so the only possibility is $\alpha(j) \equiv j \equiv -j$. By Lemma~\ref{lem:CommonNeighborjisn/2}, $N(i) \cap N(j) = \{0, i+j\}$. Applying $\alpha$ and our assumptions,
$N(\alpha(i)) \cap N(j) = \{0,\alpha(i+j)\}.$
If $\alpha(i)\equiv i$, then by Lemma~\ref{lem:CommonNeighborjisn/2} $\alpha(i+j) \equiv i+j \equiv \alpha (i) + \alpha(j)$. If instead $\alpha(i)\equiv -i$, then by Lemma~\ref{lem:CommonNeighborjisn/2},
\[
N(-i) \cap N(j) = \{0, -i+j\} = \{0, \alpha(i) + \alpha(j)\}= \{0, \alpha(i+j)\}.
\]

Now that we have shown that equation~\ref{eqn:induction} holds for $m \in \{1, 2\}$, we move on to the inductive step. Let $m \ge 3$ and assume $\alpha(ci+dj) \equiv c\alpha(i) + dj$ for all $c+d \in \{m-1, m-2\}$, with $c \ge 0$ and $d \in \{0, 1\}$. We will show it holds for $c+d=m$. 

First, suppose $d=0$. By Lemma~\ref{lem:CommonNeighborjisn/2}, $
N((m-2)i) \cap N(mi) = \{(m-1)i\}.$
Applying $\alpha$ and the inductive hypothesis, $
N((m-2)\alpha(i)) \cap N(\alpha(mi)) = \{(m-1)\alpha(i)\}.$
Since there is exactly one common neighbor, either $\alpha(mi) \equiv (m-2) \alpha(i) + 2i$ or $\alpha(mi) \equiv (m-2) \alpha(i) - 2i$. In the former case, $(m-1) \alpha(i) \equiv (m-2)\alpha(i)+i$, so $\alpha(i) \equiv i$ and thus $\alpha(mi) \equiv m \alpha(i)$. 
In the latter case, 
$(m-1) \alpha(i) \equiv (m-2)\alpha(i) -i$, so $\alpha(i) \equiv -i$, and thus $\alpha(mi) \equiv m(-i) \equiv m \alpha(i)$. Thus equation~\ref{eqn:induction} holds when $d=0$.

Next, suppose $d=1$. We must show $\alpha((m-1)i + j) \equiv (m-1)\alpha(i) + j$. By Lemma~\ref{lem:CommonNeighborjisn/2}, $N((m-2)i) \cap N((m-1)i+j) = \{ (m-1)i, (m-2)i+j\}.$
Applying $\alpha$ and using the inductive hypothesis, \begin{equation}\label{eqn:jisn/2tis1}
N((m-2)\alpha(i)) \cap N(\alpha((m-1)i+j)) = \{ (m-1)\alpha(i), (m-2)\alpha(i)+j\}.
\end{equation}
Because these two vertices have exactly two common neighbors, by Lemma~\ref{lem:CommonNeighborjisn/2}, either $\alpha((m-1)i+j)\equiv (m-2)\alpha(i) + i +j$ or $\alpha((m-1)i+j)+i+j\equiv (m-2)\alpha(i)$.

In the former case, by Lemma~\ref{lem:CommonNeighborjisn/2} the two common neighbors are $(m-2)\alpha(i) +i$ and $(m-2)\alpha(i) +j.$
Comparing this with equation~\ref{eqn:jisn/2tis1} forces $(m-2)\alpha(i) +i\equiv(m-1)\alpha(i)$ and so $\alpha(i) \equiv i$. Then $\alpha((m-1)i+j) \equiv (m-2)\alpha(i) + i +j\equiv (m-1)\alpha(i) + j.$ In the latter case, the two common neighbors are $\alpha((m-1)i+j) +i$ and $\alpha((m-1)i+j) +j$.
Adding $j$ to both sides of $\alpha((m-1)i+j)+i+j\equiv (m-2)\alpha(i)$ yields $\alpha((m-1)i+j) +i \equiv (m-2) \alpha(i) + j$.
Comparing this with~\ref{eqn:jisn/2tis1} forces
$\alpha((m-1)i+j) +j \equiv (m-1)\alpha(i)$ which means $\alpha((m-1)i+j) \equiv (m-1)\alpha(i)+j$, as desired.

This finishes the induction proof in the case $j = n/2$.

Now assume $j<n/2$. Recall we have assumed $C_n(i,j)$ is twin-free and so by Lemma~\ref{lem:whenTwins}, $i+j \neq n/2$. By Lemma~\ref{lem:CommonNeighbors}, the sets of common neighbors are affected by the presence of special conditions satisfied by $i$ and $j$. Here, we address the case where none of the special conditions hold. Some adjustments to the induction are necessary when special conditions do hold; details can be found in Appendix~\ref{sec:AutoSpecialConditions}.

We first consider a second base case, $m=2$. 
Applying $\alpha$ to $N(0) \cap N(2i) = \{i\}$ gives $N(0) \cap N(\alpha(2i)) = \{\alpha(i)\}$. Since 0 and $\alpha(2i)$ have exactly one common neighbor, by Lemma~\ref{tab:CommonNeighbors}, $\alpha(2i)$ is in $\pm 2i, \pm 2j$. Moreover, that common neighbor $\alpha(i)$ is one of $\pm i, \pm j$ and specifically such that $2 \alpha(i) \equiv \alpha(2i)$. Considering $N(0) \cap N(2j)$ similarly yields $2 \alpha(j) \equiv \alpha(2j)$. Applying $\alpha$ to $N(0) \cap N(i+j) = \{i, j\}$ gives $N(0) \cap N(\alpha(i+j)) = \{\alpha(i),\alpha(j)\}$. Since 0 and $\alpha(i+j)$ have exactly two common neighbors, by Lemma~\ref{lem:CommonNeighbors}, $\alpha(i+j)$ is one of $i+j, i-j, -i-j, -i+j$. Moreover, the only possibilities for their common neighbors $\alpha(i), \alpha(j)$ satisfy $\alpha(i) + \alpha(j) \equiv \alpha(i+j)$.

Now let $m \ge 3$ and assume $\alpha(ci+dj) \equiv c\alpha(i) + d\alpha(j)$ for all $c+d \in \{m-1, m-2\}$, with $c,d \ge 0$; we will show this holds for $c+d=m$. We have either $c \ge 2$ or $d \ge 2$. Assume $c \ge 2$.
From Lemma~\ref{lem:CommonNeighbors}, $N((c-2)i + dj) \cap N(ci+dj) = \{(c-1)i +dj\}.$ Applying $\alpha$ and the inductive hypothesis yields $N((c-2)\alpha(i) + d\alpha(j)) \cap N(\alpha(ci+dj)) = \{(c-1)\alpha(i) +d\alpha(j)\}.$ 

Lemma~\ref{lem:CommonNeighbors} shows that whenever vertices $u$ and $v$ have exactly one common neighbor $w$, the relationship $v\equiv 2w-u$ holds. 
Hence,
\[
\alpha(ci+dj) \equiv2[(c-1)\alpha(i) +d\alpha(j)] - [(c-2)\alpha(i) + d\alpha(j)] \equiv c\alpha(i) + d \alpha(j).
\]
If instead $d\ge 2$, an analogous argument on $N(ci + (d-2)j) \cap N(ci+dj) = \{ci +(d-1)j\}$ verifies equation~\ref{eqn:induction}.

This completes the proof when no special conditions hold.
\end{proof}

\begin{thm}\label{thm:TwinFreeAutG}
If $C_n(i,j)$ is not $C_{10}(1,3)$, is connected, and twin-free, then $\aut(C_n(i,j)) = \Z_n \rtimes H$, where the action of $(s, t) \in \Z_n \rtimes H$ on a vertex of $C_n(i,j)$ is $(s,t)\cdot(a) = s + ta$. 
\end{thm}

\begin{proof}

Let $\gamma \in \aut(C_n(i,j))$. First assume that $\gamma$ fixes vertex $0$. By Proposition~\ref{prop:autoFix0}, $\gamma$ induces an automorphism of the group $\Z_n$. It is well known that $\aut(\Z_n) = U(n)$. Thus, there exists $t \in U(n)$ such that $\gamma(a) = ta$. Since $\gamma$ is an automorphism of $\Z_n$, $\gamma(a) - \gamma(b) \equiv ta - tb \equiv t(a-b)$. Thus, considering $a$ and $b$ that are neighbors, we get that multiplication by $t$ must preserve $\{\pm i, \pm j\}$ and so $t \in H$. 

Next assume that $\gamma(0) = s$. Let $\sigma_{-s}$ be the translation defined by $\sigma_{-a}(a) = -s+a$. Then $\sigma_{-s} \circ \gamma$ is an automorphism of $C_n(i,j)$ that fixes $0$ and hence there is a $t \in H$ such that $\sigma_{-s} \circ \gamma(a) = ta$. But then $\gamma(a) = s+ta$ and we can represent $\gamma$ with the ordered pair $(s, t)$. 

Next we show $\aut(C_n(i,j))$ has the structure of the semidirect product $\Z_n \rtimes H$. We just showed that any $\gamma \in \aut(C_n(i,j))$ can be written as the composition of a translation by an element of $\Z_n$ after multiplication by an element of $H$. Clearly $\Z_n \cap H$ contains only the identity automorphism. We noted earlier that for all $s \in \Z_n$ and $t \in H$, $(\tau_t)^{-1} \circ \sigma_s \circ \tau_t = \sigma_{st^{-1}}$.
Thus $\Z_n$ is a normal subgroup of $\aut(C_n(i,j))$.
\end{proof}

This result aligns with Godsil's Lemma~2.2  in~\cite{G1981}, because in the twin-free case, except $C_{10}(1,3)$,  $\mathbb Z_n$ is a normal subgroup of $\aut(C_n(i,j))$, so its normalizer is the entire automorphism group.

\subsection{Automorphisms of $C_n(i,j)$ with Twins}

If $C_n(i,j)$ has twins, then by Corollary~\ref{cor:TwinsEdgeTrans}, $C_n(i,j)$ is edge-transitive and hence arc-transitive. In \cite{L2005}, Li provides a description of the automorphism group of any arc-transitive circulant graph, based on its tensor-lexicographic decomposition into a normal circulant graph, some complete graphs and an empty graph. Here we take a more elementary approach for the special case of two-generator circulant graphs. 

By Lemma~\ref{lem:whenTwins}, the only two-generator circulant graphs with adjacent twins are $C_4(1,2) \cong K_4$ and $C_5(1,2) \cong K_5$. The automorphism groups of these graphs are the symmetric groups $S_4$ and $S_5$ respectively. It can be verified computationally that the automorphism group of the one two-generator graph with co-twins, $C_{10}(1,3)$, is $\mathbb Z_2 \times S_5$.

If $C_n(i,j)$ is twin-free, then $\aut(C_n(i,j)) = \Z_n \rtimes H$ by Theorem~\ref{thm:TwinFreeAutG}. If $C_n(i,j)$ has nonadjacent twins, then $\Z_n \rtimes H$ is still a subgroup of the automorphism group. 
Additionally, $C_n(i,j)$ has automorphisms that permute 
mutual twin vertices.
Recall from Section~\ref{sec:TwinCase} that for any graph $G$, we can collapse sets of mutually twin vertices to define the twin quotient graph $\widetilde G$. Define $\pi:\aut(G) \to \aut(\widetilde G)$ by
$[\pi(\alpha)]\cdot [x] = [\alpha(x)]$
for all $\alpha \in \aut(G)$ and $x \in V(G)$. 
Properties of automorphisms guarantee that $\pi(\alpha)$ is a bijection that respects adjacency and nonadjacency in $\widetilde G$.
Note that $\ker(\pi)$ is a normal subgroup of $\aut(G)$ consisting of automorphisms that simply permute the vertices within each equivalence class.

\begin{lemma}\label{lem:sizeAutGtwins}
Let $G$ be a graph of order $n$ where every vertex is in a set of $k$ mutual twins. Then $\ker(\pi) = (S_k)^{n/k}$ and 
$
\aut(\widetilde G) = \aut(G) / \ker(\pi).
$
Hence
\[
|\aut(G)| = |\ker(\pi)| |\aut(\widetilde G)| = (k!)^{n/k}|\aut(\widetilde G)|.
\]
\end{lemma}

\begin{proof} It suffices to show that $\pi$ is surjective and then apply the First Isomorphism Theorem for groups (see~\cite{GL2018}). For each equivalence class of vertices in $G$, we select a class representative and label the vertices with subscripts in some order, $[x] = \{x_1, x_2, \dots , x_k\}$. Suppose $\widetilde \beta \in \aut(\widetilde G)$. If $\widetilde \beta([x])= [y]$, then let $\beta(x_m) = y_m$ for each $m \in \{1, 2, \dots, k\}$. 
It is easy to verify that $\beta \in \aut(G)$ and $\pi(\beta) = \widetilde \beta$. 
\end{proof} 

We can apply Lemma~\ref{lem:sizeAutGtwins} to connected, two-generator circulant graphs with twins. The relevant twin quotient graphs are $K_1$, $K_2$ and $C_{n/2}$. We find $|\aut(C_n(1,2))| = (n!)^{n/n} \cdot 1$ for $n\in \{4,5\}$, $|\aut(C_6(1,3))| = (3!)^{6/3} \cdot 2 = 72$, $|\aut(C_8(1,3))| = (4!)^{8/4} \cdot 2=1152$ and for $i+j = n/2$, but $n\neq 8$, $|\aut(C_n(i,j))| = 2^{n/2} n$.

\begin{thm}\label{thm:NonadjTwinsAutG} 
If $C_n(i,j)$ is connected and has  twins, then any automorphism of $C_n(i,j)$ is a composition of some $(s, t) \in \Z_n \rtimes H$ and an automorphism that permutes twins.
\end{thm}

\begin{proof}
When $n\in \{4,5\}$, $C_n(1,2) \cong K_n$ and the result holds. 

In general, let $N = \Z_n \rtimes H$; this is a subgroup of $\aut(C_n(i,j))$ of order $2n$. 
By what is frequently called the Second Isomorphism Theorem for groups (see~\cite{GL2018}), 
$N/(N\cap \ker(\pi)) \cong N\ker(\pi)/ \ker(\pi)$.
Comparing orders of these quotient groups, 
$
|N\ker(\pi)| = 2n |\ker(\pi)|/{|N\cap \ker(\pi)|}.
$
It is not difficult to show that 
\[
N \cap \ker(\pi) = \begin{cases}
\{(0, \pm 1), (2,\pm 1), (4, \pm 1)\}, \quad & (n, i, j) = (6, 1, 3)\\
\{(0, \pm 1), (2,\pm 1), (4, \pm 1), (6, \pm 1)\}, & (n, i, j) = (8, 1, 3)\\
\{(0,1), (n/2, 1)\}, & i+j = n/2, \text{ but } n\neq 8.
\end{cases}
\]
Simple calculations show that in each case, $|N\ker(\pi)| = |\aut(C_n(i,j))|$ and so $N\ker(\pi) = \aut(C_n(i,j))$.
\end{proof}

\section{Symmetry Parameters for Twin-Free $C_n(i,j)$}\label{sec:TwinFreeCase}

After establishing $\aut(C_n(i,j)) = \Z_n \rtimes H$ in the twin-free case, we can find the symmetry parameters with relative ease. The result below proves in the affirmative a conjecture on the determining number of connected, twin-free, two-generator circulant graphs of Brooks et al.\cite{Brooks2021}.

\begin{thm}\label{thm:DetDistCostCaseI}
If $C_n(i,j)$ is not $C_{10}(1,3)$, is connected and twin-free, then
\[
\det(C_n(i,j))=2, \, \dist(C_n(i,j)) = 2, \text{ and } \rho(C_n(i,j)) =3.
\]
\end{thm}

\begin{proof}
For determining, first let $a \in \Z_n$. Then $a$ is fixed by the nontrivial automorphism $(2a, -1)$, so $\{a\}$ is not a determining set. Thus $\det(C_n(i,j)) > 1.$ Next, let $W = \{0, 1\}$ and assume $\alpha = (s, t) \in \aut(C_n(i,j))$ fixes both vertices in $W$.
Then $0 \equiv \alpha(0) \equiv s+t\cdot 0 \equiv s$. Next $1 \equiv \alpha(1) \equiv s +t \cdot 1 \equiv 0+t \equiv t$. Hence $\alpha= (0, 1)$ which is the identity and $W$ is a determining set. Thus $\det(C_n(i,j))= 2$.

Next we find a $2$-coloring that is distinguishing and that has one color class of size $3$. There are two cases to consider.

\medskip

\noindent {\bf Case 1.} Assume vertices $i$ and $j$ are not adjacent in $C_n(i,j)$. If $C_n(i,j)$ is not edge-transitive, color the vertices in $\{0, i, j\}$ red and all other vertices blue. Assume $\alpha = (s,t) \in \aut(C_n(i,j))$ is an automorphism that preserves the color classes. Since $\{0,i,j\}$ induces a path, $\{\alpha(i), \alpha(j)\} = \{i, j\}$ and $\alpha(0) = 0$.
Hence $0 \equiv s +0t \equiv s$. 
For $\alpha$ to be nontrivial, $t \not \equiv 1$ and thus $\alpha(i)\equiv ti \equiv j$ and $\alpha(j)\equiv tj \equiv i$. 
Since $t \in H$, $ti \equiv \pm i$. 
Thus $j \equiv \pm i$, or equivalently either $i+j \equiv 0$ or $j-i \equiv 0$. Our assumption throughout this paper is that $0<i<j\le n/2$, which implies both $0<i+j < n$ and $0<j-i <n$, a contradiction. Hence $\alpha$ must be the identity and so $\{0, i, j\}$ is a color class in a 2-distinguishing coloring.

In the case $C_n(1,j)$, $j^2 \equiv \pm 1$, we have $H = \{\pm 1, \pm j\}$. When $j^2 \equiv -1$, $\{0,i,j\}$ is still a color class in a 2-distinguishing coloring. However, when $j^2 \equiv 1$, we have that $t=j$ gives a nontrivial automorphism that preserves this set. In this case, by a similar argument, $\{0,-i,j\}$ is instead a color class in a 2-distinguishing set.

\medskip

\noindent {\bf Case 2.} If $i$ and $j$ are adjacent in $C_n(i,j)$, then $i \in N(j) = \{0, i+j, -i+j, 2j\}$. 
The assumption on $i$ and $j$ imply $0 < i <i+j < 2j \le n$, so $i \equiv -i+j $. Equivalently, $2i \equiv j$. The only edge-transitive two-generator circulant graphs satisfying this condition are $C_5(1,2)$ and $C_6(1,2)$, both of which have twins. Hence in this case, $C_n(i, j)$ is not edge-transitive.

In this case, we can show 
$ -j \notin N(i) = \{0, 2i, i+j, i-j\} = \{0, j, i+j, i-j\}$.
Since $0<j \le n/2$, $j \not \equiv 0$, so $-j \not \equiv 0$.
 If $-j \equiv 2i \equiv j$, then $j = n/2$, so $4i \equiv 0$ and $\gcd(n, i, j) = i$. Since we assumed $C_n(i,j)$ is connected, $\gcd(n, i, j) = 1$, which forces $n = 4$, contradicting our assumption that $C_n(i,j)$ is twin-free. 
 If $-j \equiv i+j$, then $0 \equiv i+2j \equiv 5i$ and again we get $\gcd(n, i, j) = i$. Thus $i=1$ and $n=5$, again contradicting the assumption that $C_n(i,j)$ is twin-free. 
 Finally, if $-j \equiv i-j$, then $i \equiv 0$, contradicting $0<i<j \le n/2$. Thus, $-j$ is not adjacent to $i$.

We now color the vertices in $\{-j, 0, i\}$ red and all other vertices blue. Arguing as above, we find that any automorphism preserving the color classes must be trivial.

\medskip

In both cases, $\dist(C_n(i,j)) = 2$ and $\rho(C_n(i,j)) \le 3$.
A color class in any 2-distinguishing coloring cannot be a singleton set because $\det(C_n(i,j)) = 2$. 
Furthermore, if $a\neq b$ in $\Z_n$, then the nontrivial automorphism $(a+b, -1)$ interchanges them, so a color class in a 2-distinguishing coloring cannot consist of just two vertices. Thus $\rho(C_n(i,j)) = 3$.
\end{proof}

\section{Subdivided Circulant Graphs}\label{sec:SubdividedCase}

We next consider the symmetry parameters of the subdivided circulant graphs $C_n(i_\div p,j)$ and $C_n(i, j_\div p)$.
Recall that these are the derived graphs associated to a bouquet voltage graph $B_2$ in which one of the arcs is subdivided by $p>1$ vertices of degree $2$. 
Because the voltage graph has order at least $2$, we can no longer label vertices of the derived graph simply with elements of $\Z_n$. As shown in Figure~\ref{fig:SubdividedVoltage}, we label them $u_a$ and $v_a^r$ 
where $a \in \Z_n$ and $r \in \{1, \dots, p\}$.

\subsection{$C_n(i_\div p, j)$ and $C_n(i, j_\div p)$, $j<n/2$} \label{subsec:jlessn/2}

We now consider the case in which the arc with voltage $i$ has been subdivided by $p$ vertices. In the derived graph, $u_a$ is no longer adjacent to the (distinct) vertices $u_{a+i}$ and $u_{a-i}$.
Instead, for each $a \in \Z_n$, there is a path $
(u_a, v_a^1, v_a^2, \dots, v_a^p, u_{a+i}).
$ of length $p+1$, all of whose interior vertices have degree 2. 
Note that there is no ambiguity regarding the subscripts on the degree-2 vertices: if $b \equiv a+i$, it cannot also be the case that $a \equiv b+i$, because $0< 2i < n$.
By definition,
\[
N(u_a) = \begin{cases}
\{v_a^1, v_{a-i}^p, u_{a-j}, u_{a+j}\}, \quad & j<n/2,\\
\{v_a^1, v_{a-i}^p, u_{a+j}\}, & j=n/2.
\end{cases}
\]
Thus each $u_a$ has a distinct pair of degree-$2$ neighbors, meaning that no two vertices of this type are twins. 
Additionally, each $v_a^r$, $r \in \{1, \dots p\}$, is uniquely identified by its distances from $u_a$ and $u_{a+i}$. Hence $C_n(i_\div p, j)$ is twin-free.
Neighborhood sizes dictate that the only candidates for subdivided graphs with co-twins are small and direct inspection shows none exist.

Clearly, $C_n(i_\div \ell, j)$ is neither vertex-transitive nor edge-transitive. However, automorphisms of $C_n(i,j)$ extend uniquely to automorphisms of $C_n(i_\div \ell, j)$, provided that they do not interchange edges of the form $\{a, a+i\}$ and $\{b, b+j\}$. We let
\[
H' = H \cap \aut(\mathbb Z_n, \{\pm i\},\{\pm j\}) = 
\begin{cases} 
\{\pm 1\}, \quad & C_n(i,j) \text { is edge-transitive},\\
H, & \text{otherwise.}
\end{cases}
\]

\begin{lemma}\label{lem:aut_i_extension}
 For any $\alpha' = (s,t)  \in \mathbb Z_n \rtimes H'$,
 there is a unique $\alpha \in \aut(C_n(i_\div \ell,j))$ such that $\alpha(u_a) = \alpha'(u_a)$ for all $a \in \mathbb{Z}_n$. The action of the unique extension $\alpha$ is defined as: 
 \begin{equation}\label{eqn:AutExti}
 \alpha(x)= \begin{cases}
 u_{s+ta}, \quad & x = u_a, \\ 
 v_{s +at}^r, \quad &x = v_a^r \text{ and } ti\equiv i, \\
 v_{s+ta-i}^{(\ell+1)-r}, &x = v_a^r \text{ and } ti \equiv -i.
 \end{cases}
 \end{equation}
\end{lemma}

\begin{proof} 
To show uniqueness, suppose $\gamma, \lambda\in \aut(C_n(i_\div p, j))$ satisfy $\gamma(u_a) = \lambda(u_a)$ for all $a \in \Z_n$. Then $\gamma^{-1} \circ \lambda$ is an automorphism of $C_n(i_\div p, j)$ that fixes all of the non-degree-$2$ vertices. Since each $v_a^r$ is uniquely identified by its distances from $u_a$ and $u_{a+i}$, all the degree-$2$ vertices must also be fixed by $\gamma^{-1} \circ \lambda$. Thus $\gamma^{-1} \circ \lambda$ is the identity, and so $\gamma = \lambda$.

To define that action of $\alpha$ on the degree-2 vertices, we consider cases based $ti \equiv i$ or $ti \equiv -i$. In the former case, $(s, t)\cdot (u_{a+i}) = u_{s+t(a+i)} = u_{(s+ta) +i}$. In this case, we map the path of degree-2 vertices between $u_a$ and $u_{a+i}$ to the path of degree-2 vertices between $u_{s+ta}$ and $u_{(s+ta)+i}$, in the same order. In the latter case, $(s, t)\cdot (u_{a+i}) = u_{s+(a+i)t} = u_{(s+ta) -i}$. In this case, we map the path of degree-2 vertices between $u_a$ and $u_{a+i}$ to the path of degree-2 vertices between $u_{s+ta}$ and $u_{(s+ta)-i}$, in `reverse order.' It can be easily verified that $\alpha$ is respects adjacency and nonadjacency by checking that the action on the three types of edges: $u_a u_{a\pm j}$ for $a \in \mathbb{Z}_n$, $v_{a}^r v_a^{r\pm 1}$ for $a \in \mathbb{Z}_n$, and $u_{a} v_{b}^r$ with $a \in \mathbb{Z}_n$ and $(b,r) \in \{(a,1), (a-i,p)\}$.
\end{proof}

 \begin{thm}\label{thm:AutII} 
Assume $C_n(i,j)$ is connected. Then
$\aut(C_n(i_\div p, j))=\Z_n \rtimes H'$, with the action of elements of $\Z_n \rtimes H'$ as defined in Equation~\ref{eqn:AutExti}.
 \end{thm}

 \begin{proof}
By Lemma~\ref{lem:aut_i_extension}, we have $\Z_n \rtimes H' \subseteq \aut(C_n(i_\div p, j))$. Thus, we show that $\aut(C_n(i_\div p, j))\subseteq \Z_n \rtimes H'$. Let $\alpha \in \aut(C_n(i_\div p, j))$. Since automorphisms respect degree, $\alpha$ restricts to a bijection $\alpha'$ on the set of non-degree-$2$ vertices, $\{u_a \mid a \in \Z_n\}$. By definition, $u_a$ and $u_b$ are adjacent as vertices in $C_n(i,j)$ if and only if, as vertices in $C_n(i_\div p,j)$, they are either adjacent or there is a unique path between them of length $p+1$, all of whose interior vertices have degree $2$. These are properties respected by the automorphism $\alpha$.
 Thus $u_a$ and $u_b$ are adjacent as vertices in $C_n(i,j)$ if and only if $\alpha'(u_a)$ and $\alpha'(u_b)$ are adjacent as vertices in $C_n(i,j)$. Hence $\alpha'$ is an automorphism of $C_n(i,j)$.

 If $C_n(i,j)$ is twin-free and not $C_{10}(1,3)$, then by Theorem~\ref{thm:TwinFreeAutG}, $\alpha' \in \aut(C_n(i,j)) = \Z_n \rtimes H'$ and we are done. If $C_n(i,j)$ has twins, then by Theorem~\ref{thm:NonadjTwinsAutG}, $\alpha'$ is
 the composition of an element of $\mathbb Z_n \rtimes H'$ and an automorphism $\rho$ that permutes twins. In each of the four special cases of Lemma~\ref{lem:whenTwins}, namely $C_4(1,2) = K_4$, $C_5(1, 2) = K_5$, $C_6(1,3)$ and $C_8(1,3)$, 
 as well as in the co-twin case $C_{10}(1,3)$,
 the edges corresponding to arcs of voltage $i=1$ induce a Hamilton $n$-cycle in $C_n(i,j)$. Moreover, $\alpha'$ respects this $n$-cycle because $\alpha$ respects the Hamilton cycle in $C_n(i_\div p, j))$ induced by all edges having at least one endvertex of degree 2. Thus $\alpha'$ is an automorphism of an $n$-cycle and so $\rho$ is trivial.

 Finally, suppose $C_n(i,j)$ has twins because $i+j = n/2$, but $n \neq 8$. In this case, each vertex $a$ has exactly one nonadjacent twin, namely $a+i+j$. For each $a$, let $\chi_a$ be the automorphism of $C_n(i,j)$ that exchanges $a$ and $a+i+j$ and leaves all other vertices fixed. Note that $\chi_a = \chi_{a+i+j}$.
 By Theorem~\ref{thm:NonadjTwinsAutG}, $\alpha'=(s,t)\circ \rho$ where $(s,t) \in \Z_n \rtimes H'$ and $\rho = \chi_{a_1} \circ \cdots \circ \chi_{a_k}$.
 We will show that if $\rho$ exchanges one pair of twins, then $\rho$ must exchange all pairs of twins, and so $\rho = (n/2, 1) \in \Z_n \rtimes H'$. We begin by noting that both $\alpha$ and $(s,t)$ are elements of $\aut(C_n(i_\div p, j))$, and hence so is $(-s, t^{-1})\circ \alpha.$ By construction, 
 \[[(-s, t^{-1})\circ \alpha] \cdot (u_a) = (-s, t^{-1})\cdot(\alpha'(u_a)) = (-s, t^{-1})\cdot((s,t)\cdot \rho(u_a)) = \rho(u_a) .\]
 If $\chi_a$ is a factor of $\rho$, then $\rho(u_a)= u_{a+i+j}$. Then $\rho$ must take the degree-4 neighbors of $u_a$ to the degree-4 neighbors of $u_{a+i+j}$. This means
 $\{\rho(u_{a+j}), \rho(u_{a-j})\} = \{u_{a+i+2j}, u_{a+i}\}$.
 Hence, both $\chi_{a+j}$ and $\chi_{a+i}$ must be factors of $\rho$ as well. Repeating this argument, we find that for all positive integers $x$ and $y$, $\chi_{a+xi+yj}$ is a factor of $\rho$. Since our assumptions imply $\gcd(n,i,j) = 1$, we thus have $\rho = \chi_0 \circ \cdots \circ \chi_{\frac{n}{2}-1}$. 
 \end{proof}

 \begin{prop}\label{prop:IInospecialt} 
 If $H'= \{\pm 1\}$ then a minimum determining set of $C_n(i_\div p, j)$ is $W=\{v_0^1\} $ if $p \ge 2$ and $W= \{u_0, v_0^1\}$ if $p=1$.
 \end{prop}
 
 \begin{proof}

 First assume $p\ge 2$. Suppose $(s, t)\in \aut(C_n(i_\div p,j))$ fixes $v_0^1$. We must show $(s,t)$ is the identity; that is, $s\equiv 0$ and $t\equiv 1$ in $\Z_n$.
 If $t \not \equiv 1$, then $t \equiv -1$ so by Equation~\ref{eqn:AutExti} we have $(s,t) \cdot (v_0^1)=v_{s-i}^p$. Since $v_0^1$ is fixed, this implies $p=1$, contradicting our assumption. Hence $t \equiv 1$, and $(s,t) \cdot (v_0^1) = v_s^1$ yields that $s \equiv 0$.

 Next assume $p=1$. 
Suppose $(s,t)$ fixes both $v_0^1$ and $u_0$. Since $(s,t)$ fixes $u_0$, $s\equiv 0$. 
 If $t \not \equiv 1$, then $t \equiv 1$ and by Equation~\ref{eqn:AutExti}, we have $(s,t)\cdot (v_0^1)=v_{-i}^1.$
 Since $0<i<n/2$, $-i \not \equiv 0$, contradicting our assumption that $(s,t)$ fixes $v_0^1$. Hence $t\equiv 1$ and $\{u_0, v_0^1\}$ is a determining set. 

 For minimality, any determining set is nonempty as $C_n(i_\div p, j)$ has nontrivial automorphisms. In the case $p=1$, neither  $\{u_a\}$ nor  $\{v_a^1\}$ can be determining as they are fixed by the nontrivial automorphisms $(2a, -1)$ and $(i, -1)$ respectively.
 \end{proof}

Now assume $H' \neq \{\pm 1 \}$. This means $C_n(i,j)$ is not edge-transitive and  $H'= \{ \pm 1, \pm t\}$ for some $1 \neq t \in U(n)$ such that $ti \equiv i$ and $tj \equiv -j$. Also, neither $i$ nor $j$ is a unit. We require some additional algebraic results.

\begin{lemma}\label{lem:units_nonunits}
 Let $0 < i < n/2$. If $i \not \in U(n)$, then there exists a prime dividing $n$ that does not divide $i$ if and only if there exists $b$ such that $b \not \in U(n)$ but $b-i \in U(n)$. 
\end{lemma}

\begin{proof} 
First assume there exists at least one prime dividing $n$ that does not divide $i$. Let $b$ be the product of all such primes. Then $\gcd(n, b) = b > 1$, so $b \not \in U(n)$. 
Since no prime dividing $n$ divides both $b$ and $i$, $\gcd(n, b-i) = 1$, meaning $b-i \in U(n)$.

Conversely, assume that every prime dividing $n$ also divides $i$. For any $b \not \in U(n)$, there is some prime dividing both $n$ and $b$. By assumption, this prime also divides $i$, and thus divides $b-i$. Thus $b-i \not \in U(n)$.
\end{proof}

 \begin{cor}\label{cor:unit_nonunit}
 Let $n, i, j$ satisfy $\gcd(n, i, j) = 1$ and $0 < i < j \le n/2$. If $H' \neq \{\pm 1\}$, then there exists $a \in U(n)$ such that $a+i \notin U(n)$.
\end{cor}

\begin{proof}
 Assume $H' \neq \{\pm 1\}$. 
Since $j$ is not a unit, there is a prime  dividing $n$ that also divides $j$. However, this prime cannot divide $i$ because $\gcd(n, i, j) = 1$. Hence by Lemma~\ref{lem:units_nonunits}, there exists $b \not \in U(n)$ such that $b-i \in U(n)$. Let $a \equiv b-i$.
\end{proof}

\begin{prop}\label{prop:IIspecialt}
 Let $n, i, j$ satisfy $\gcd(n, i, j) = 1$ and $0 < i < j \le n/2$, and assume $H' \neq \{\pm 1\}$.
 Let $a \in U(n)$ such that $a+i \notin U(n)$. Then $W= \{u_0, v_a^1\}$ is a minimum determining set of $C_n(i_\div p,j)$. 
\end{prop}

\begin{proof}
Suppose $(s,t) \in \aut(C_n(i_\div p,j))$ fixes $\{u_0, v_a^1\}$. Since $(s,t)$ fixes $u_0$, $s\equiv 0$. By Equation~\ref{eqn:AutExti}, $(0,t)\cdot (v_a^1)$ is $v_{at}^1$ or $v_{at-i}^{p}$ depending on whether $ti \equiv i$ or $ti \equiv -i$. In the former case, since $v_a^1$ is fixed, $a\equiv ta$, so $a(1-t) \equiv 0$. By assumption, $a$ has a multiplicative inverse, so this implies $1-t \equiv 0$. Hence $(s,t) = (0,1)$. In the latter case, $p=1$ and $ta -i \equiv a$. By substitution, $ta +ti \equiv a$, so $t(a+i) = a$. By assumption, $a, t \in U(n)$, but $a+i \not \in U(n)$, a contradiction. 
 Thus the only automorphism fixing $W= \{u_0, v_a^1\}$ is the identity automorphism $(0, 1)$.

 We have already seen that for any $a \in \Z_n$, $(2a, -1)$ is a nontrivial automorphism fixing $u_a$, so 
$\{u_a\}$ is not a determining set.
By Corollary~\ref{cor:Hnot1-1}, there exists $t' \in H'$ such that $t'\not \equiv 1$ but $t'i \equiv i$.
For any $a \in \Z_n$, $(a-at',t')$ is a nontrivial automorphism fixing $v_a^r$ for all $r \in \{1, 2, \dots , p\}$, and so 
$\{v_a^r\}$ is not a determining set. 
Thus $\det(C_n(i_\div p,j) \ge 2$.
 \end{proof}

\begin{thm}\label{thm:CaseII} 
 If $p \ge 2$ and $H'=\{\pm 1\}$, then
\[
\det(C_n(i_\div p,j)) =1, \,\dist(C_n(i_\div p,j)) = 2 \text{ and }\rho(C_n(i_\div p,j))=1.
\]
Otherwise,
$
\det(C_n(i_\div p,j))= \dist(C_n(i_\div p,j)) = \rho(C_n(i_\div p,j)) = 2.
$
\end{thm}

\begin{proof}
 If $p \ge 2$ and $H'= \{\pm 1\}$, then by Proposition~\ref{prop:IInospecialt}, $C_n(i_\div p,j)$ has a one-element determining set. Any graph with determining number $1$ has distinguishing number $2$ and cost $1$, as shown in~\cite{BCKLPR2020b}. 
 
 If $p=1$ and $H'=\{\pm 1\}$, then by Proposition~\ref{prop:IInospecialt}, $W= \{u_0, v_0^1\}$ is a minimum determining set. If $H' \neq \{\pm 1\}$, then by Proposition~\ref{prop:IIspecialt}, there exists $a \in \Z$ such that $\{u_a, v_0^1\}$ is a minimum determining set. In each case, the two vertices in the determining set have different degrees. Thus, coloring these vertices red and the other vertices blue produces a 2-distinguishing coloring with cost 2.
 Since $\det(C_n(i_\div p,j))=2$, no distinguishing $2$-coloring can have a color class of size $1$.
\end{proof}

Our discussion of $C_n(i_\div p, j)$ relies several times on the fact that $i \neq n/2$ under the overall assumption $0<i<j\le n/2$. For example, this allowed us to conclude that $u_{a+i}$ and $u_{a-i}$ are distinct vertices, and that the subscripts on the degree-2 vertices $v_a^r$ are unambiguous. The other overall assumption of this paper is that $\gcd(n, i, j) = 1$, in which $i$ and $j$ play interchangeable roles. Thus our results on $C_n(i_\div p, j)$ carry over to the case $C_n(i, j_\div p)$ when $j < n/2$.

 \begin{thm} 
Assume $C_n(i,j)$ is connected and $j < n/2$. For any $\alpha' \in \aut(C_n(i,j))$, there is a unique $\alpha \in \aut(C_n(i,j_\div p))$ such that $\alpha(u_a) = \alpha'(u_a)$ for all $a \in \Z_n$. If $\alpha' = (s, t) \in \Z_n \rtimes H'$, define 
 \begin{equation}\label{eqn:AutExtj}
 \alpha(x)= \begin{cases}
 u_{s+ta}, \quad & x = u_a, \\ 
 v_{s +at}^r, \quad &x = v_a^r \text{ and } tj\equiv j, \\
 v_{s+ta-j}^{(p+1)-r}, &x = v_a^r \text{ and } tj \equiv -j.
 \end{cases}
 \end{equation}
Then
$\aut(C_n(i_\div p, j))=\Z_n \rtimes H'$, with the action of elements of $\Z_n \rtimes H'$ as defined in Equation~\ref{eqn:AutExtj}.
 \end{thm} 
 
 \begin{proof}
 In general, we can simply modify the proofs of Lemma~\ref{lem:aut_i_extension} and Theorem~\ref{thm:AutII} by interchanging $i$ and $j$. We must be more careful in the portion of the proof dealing with the case where $C_n(i,j)$ has twins
 or co-twins.
 The special cases $C_4(1,2)$ and $C_6(1,3)$ do not satisfy $j<n/2$. The other 
 special cases, $C_5(1,2)$, $C_8(1,3)$ 
and $C_{10}(1,3)$ have 
$\gcd(n,j) = 1$. 
 Thus, the edges in $C_n(i, j_\div p)$ having at least one endvertex of degree $2$ again induce a Hamilton cycle. Since any $\alpha\in \aut(C_n(i, j_\div p))$ must respect this Hamilton cycle, the restriction $\alpha'$ must respect the $n$-cycle in $C_n(i,j)$ induced by the edges corresponding to the arcs of voltage $j$, and so we can again conclude that $\alpha'$ is an element of the dihedral group $D_{2n}$.
 \end{proof}
 
 If $j < n/2$, then Proposition~\ref{prop:IInospecialt}, Lemma~\ref{lem:units_nonunits}, Corollary~\ref{cor:unit_nonunit} and Proposition~\ref{prop:IIspecialt} all hold with the roles of $i$ and $j$ interchanged. We state the main result on the values of the symmetry parameters in the case $j < n/2$ below. 
 
 \begin{thm}\label{thm:CaseIII} 
Assume $j< n/2$. If $p \ge 2$ and $H'=\{\pm 1\}$, then
$$
\det(C_n(i,j_\div p)) =1, \,\dist(C_n(i,j_\div p)) = 2 \text{ and }\rho(C_n(i,j_\div p))=1.$$
Otherwise,
$
\det(C_n(i,j_\div p))= \dist(C_n(i,j_\div p)) = \rho(C_n(i,j_\div p)) = 2.
$
\end{thm}

\subsection{$C_n(i, j_\div p)$, $j=n/2$}\label{subsec:jisn/2}

We recall some results for $C_n(i,j)$ with $j=n/2$. In this case, $n = 2j$ and $j \equiv -j$.
From Lemma~\ref{lem:whenTwins}, $C_{2j}(i,j)$ is twin-free except in two cases: $C_4(1,2)$, in which any two vertices are adjacent twins, and $C_6(1,3)$, in which $u_a$ and $u_{a+2}$ are nonadjacent twins for all $a \in \Z_6$. From Lemma`\ref{lem:CommonNeighborjisn/2}, if $C_{2j}(i,j)$ is twin-free, then for all $a \in \Z_{2j}$, we have $N(u_a) \cap N(u_{a+2i}) =\{u_{a+i}\}$ and $ N(u_a) \cap N(u_{a+i+j}) = \{u_{a+i}, u_{a+j}\}$,
and these are the only possibilities for two vertices to have common neighbors.

\begin{figure}
\includegraphics[width= 0.6\textwidth, center]{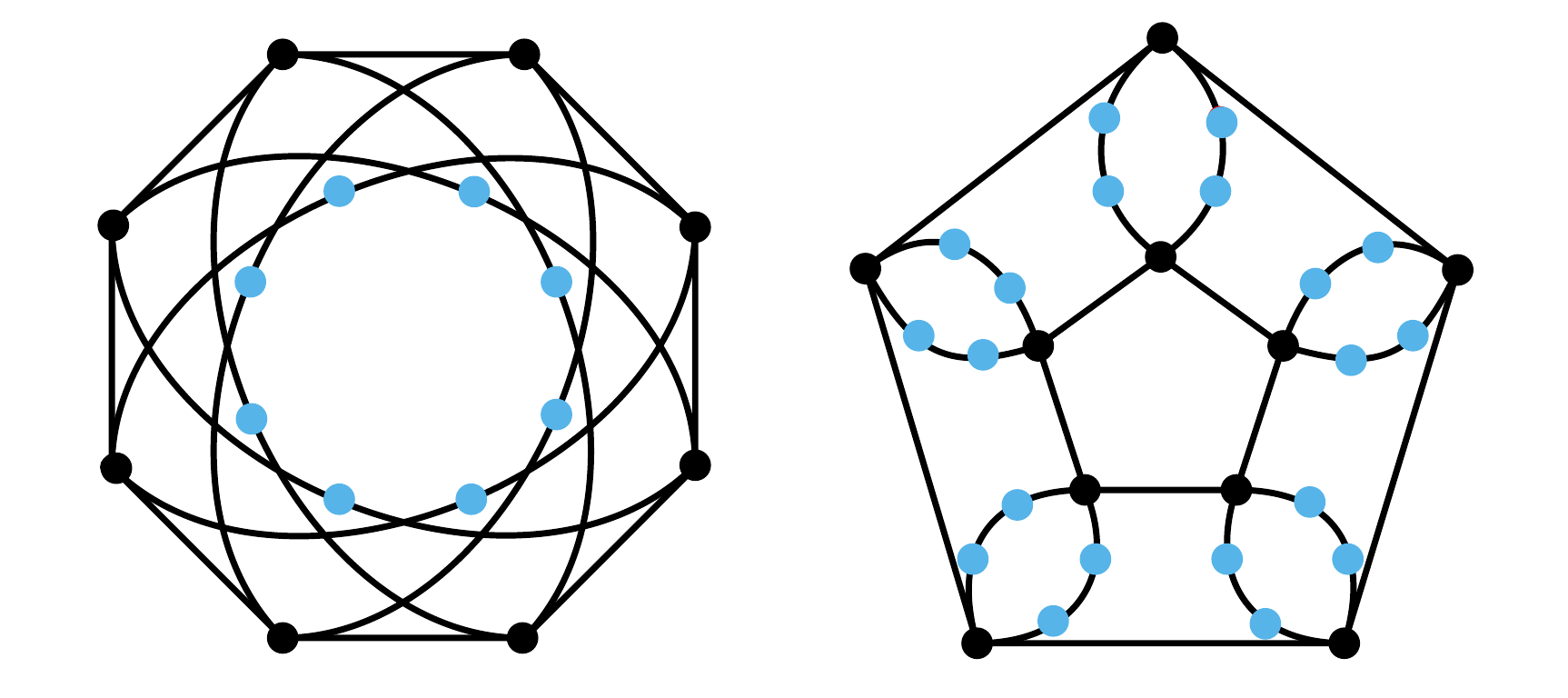}
\caption{$C_8(1,4_\div 1)$ and $C_{10}(2,5_\div 2)$.}
\label{fig:ExCaseIII}
\end{figure}

If we subdivide the loop with voltage $j$ in the voltage graph with $p$ vertices of degree 2, then in the derived graph, the single (undirected) edge between $u_a$ and $u_{a+j}$ is replaced with two paths,
\[
(u_a, v_a^1, v_a^2, \dots, v_a^p, u_{a+j}) \text{ and } (u_{a+j}, v_{a+j}^1, v_{a+j}^2, \dots , v_{a+j}^p, u_{a+2j} = u_a).
\]
Thus $N(u_a) = \{v_a^1, v_{a+j}^p, u_{a-i}, u_{a+i}\}$.
Examples $C_8(1,4_\div 1)$ and $C_{10}(2,5_\div 2)$ appear in Figure~\ref{fig:ExCaseIII}.

\begin{lemma}\label{lem:whenTwinsIII}
\hfill
\begin{enumerate}[(1)]
\item If $p \ge 2$, then $C_{2j}(i,j_\div p)$ is twin-free. 
\item If $p=1$, then for distinct $a,b$, $N(v_a^1) = N(v_b^1)$ if and only if $b = a+j$, and $N(u_a) = N(u_b)$ if and only if $(2j,i,j) = (4,1,2)$ and $b = a+2$.

\end{enumerate}
\end{lemma}

\begin{proof} 
 First assume $p\ge 2$. Then each $u_a$ has a distinct pair of degree-2 neighbors, and so no two such vertices can be twins. Moreover, each vertex $v_a^k$ has at least one neighbor of the form $v_a^{k\pm 1}$, meaning that for $a\not \equiv b$, vertices $v_a^k$ and $v_b^\ell$ cannot be twins. Finally, two vertices of the form $v_a^k$ and $v_a^\ell$ with $k < \ell$ cannot be twins because either $u_a$ (when $k=1$) or $v_a^{k-1}$ (when $k>1$) is in $ N(v_a^k)\backslash N(v_a^\ell)$.

Next assume $p=1$. Then $N(v_a^1) = \{u_a, u_{a+j}\} = N(v_{a+j}^1)$,
so $v_a^1$ and $v_{a+j}^1$ are nonadjacent twins. 
Next, assume $N(u_a) = N(u_b)$. For $u_a$ and $u_b$ to have the same degree-2 neighbors, $b = a+j$. For them to have the same degree-4 neighbors, $\{u_{a-i}, u_{a+i}\} = \{u_{a+j-i}, u_{a+j+i}\}$. Since $j \not\equiv 0$, $a+j + i \not \equiv a + i$ and $a+j- i \not \equiv a - i$. However, it is possible that $a+ i \equiv a +j - i$, namely when $2i+j \equiv 0$, in which case it also follows that $a- i \equiv a +j + i$. Under the general assumption that $\gcd(n, i, j) = 1$, $2i+j \equiv 0$ if and only if $(2j,i,j) = (4,1,2)$.
 \end{proof}
 
 \begin{figure}[h]
\includegraphics[width= 0.65\textwidth, center]{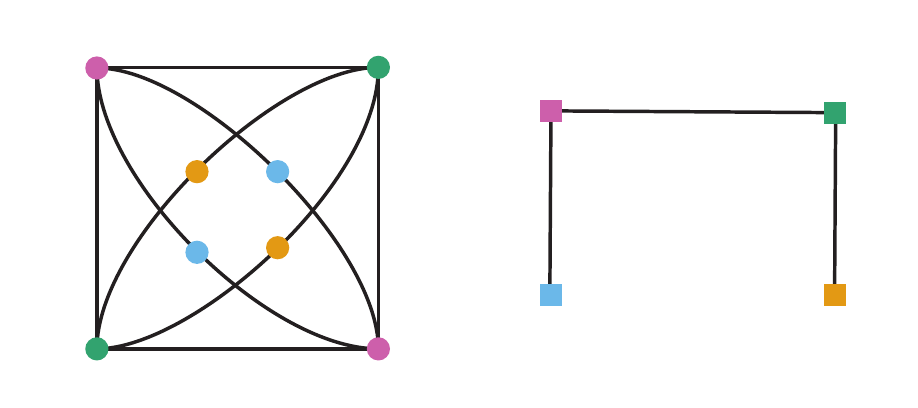}
\caption{$C_4(1,2_\div 1)$ and its twin quotient graph, $P_4$.}
\label{fig:specialCaseIII}
\end{figure}

\begin{Ex} \label{ex:C12div1}
In the special case $C_4(1,2_\div 1)$, each vertex is in a set of $k=2$ mutual nonadjacent twins and so we can use the techniques of Section~\ref{sec:TwinCase}. See Figure~\ref{fig:specialCaseIII}.
The minimum twin cover $T=\{u_0, u_1, v_0^1, v_1^1\}$ is also a minimum determining set, which implies that $\det(C_4(1,2_\div 1))=4$. The twin quotient graph is $P_4$, which has distinguishing number 2. By Theorem~\ref{thm:DistTwins}, $\dist(C_4(1,2_\div 1)) = 3$.
\end{Ex}

In general, if $p=1$ and $n = 2j\ge 6$, then the set $T = \{v_0^1, v_1^1, \dots , v_{j-1}^1\}$ is a minimum twin cover, but not a determining set. To see this, we must investigate the automorphism group of $C_{2j}(i, j_\div p)$. 

First, recall from Lemma~\ref{lem:His1-1} that when $j = n/2$, $H'= \{\pm 1\}$, so $\Z_{2j} \rtimes H'$ is the dihedral group $D_{4j}$. For all $C_{2j}(i,j)$ where $n = 2j\ge 6$ except $(n,i,j) = (6,1,3)$, $C_{2j}(i,j)$ is twin-free by Lemma~\ref{lem:whenTwins}. Then $\aut(C_{2j}(i,j)) = \Z_{2j} \rtimes H' = D_{4j}$ by Theorem~\ref{thm:TwinFreeAutG}. In the excepted case, by Theorem~\ref{thm:NonadjTwinsAutG}, every element of $\aut(C_6(1,3))$ is the composition of an element of $\Z_6 \rtimes \{\pm 1\} = D_{12}$ and an automorphism that interchanges vertices in some subset of pairs of nonadjacent twins.

When $n=2j$, $tj \equiv j$ for any $t \in H'$, so for any $a \in \Z_{2j}$ and $(s,t) \in \Z_{2j} \rtimes H'$, $(s,t)\cdot (u_{a+j}) = u_{s + ta + j}$. 
Thus when extending $(s,t)\in \Z_{2j} \rtimes H'$ to the degree-2 vertices in $C_n(i, j_\div p)$, we do not have to worry about whether to reverse the order of the superscripts, as we did when extending $(s, t)$ to the degree-2 vertices in $C_{2j}(i_\div p, j)$. That is, we can simply set $(s,t) \cdot (v_a^r) = v^r_{s+ta}$. 

\begin{figure}[h]
\includegraphics[width=0.6
\textwidth, center]{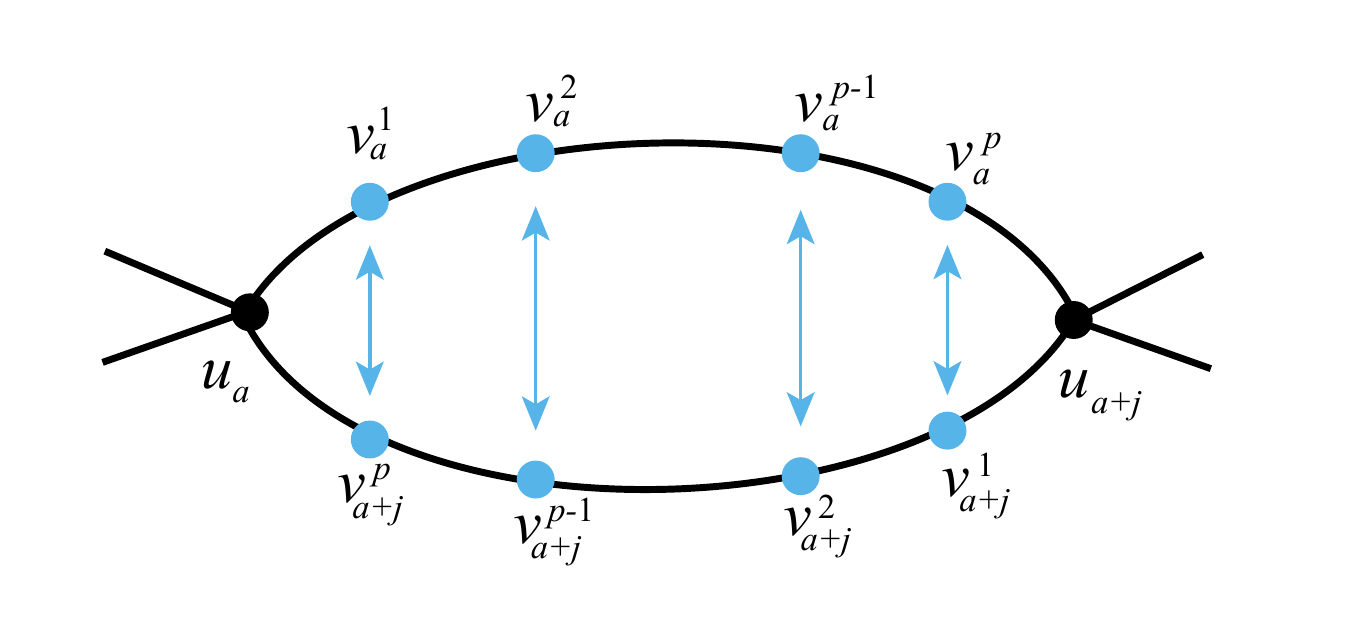}
\caption{Action of the automorphism $\beta_a = \beta_{a+j}$.}
\label{fig:miniflip}
\end{figure}

Unlike the situation for $C_{2j}(i_\div p, j)$, there are additional automorphisms of $C_{2j}(i, j_\div p)$ beyond these extensions of the automorphisms in $\Z_{2j} \rtimes H'$. For each $a \in \Z_{2j}$, there is an automorphism $\beta_a$ that `flips' the degree-2 vertices on the two paths between $u_a$ and $u_{a+j}$ and leaves all other vertices fixed; see Figure~\ref{fig:miniflip}. 
More precisely, $\beta_a \cdot (u_b) = u_b \text{ for all } b \in \Z_n$
and
\[
\beta_a \cdot (v_b^r)= \begin{cases} 
v_{b+j}^{(p+1)-r},\quad & \text{ if } b\in \{a, a+j\},\\
v_b^r, & \text{ if } b\in \Z_{2j} \backslash \{a, a+j\}.
\end{cases}
\]
Note that $\beta_a = \beta_{a+j}$, so there are exactly $j$ such automorphisms. Also, for all $a, a' \in \Z_n$, 
$\beta_a \circ \beta_{a'} = \beta_{a'} \circ \beta _a$ and $\beta_a \circ \beta_a = \iota_B$,
where $\iota_B$ is the identity. 
The subgroup generated by these automorphisms is $B= \langle \beta_a \mid 0 \le a < j\rangle \cong (\Z_2)^j$. An element $\beta \in B$ is of the form 
$\beta = \beta_0^{e_0} \cdots \beta_{j-1}^{e_{j-1}}$,
where $e_0, \dots , e_{j-1} \in \{0, 1\}$.

\begin{thm}\label{thm:AutGjisn2}
The automorphism group of $C_{2j}(i, j_\div p)$ is
\[\aut(C_{2j}(i, j_\div p)) = B \rtimes (\Z_{2j} \rtimes H') = (\Z_2)^j \rtimes (\Z_{2j} \rtimes \{\pm 1\}).\] 
An element is denoted $\beta \circ(s,t)$, where $\beta \in B$ and $(s,t) \in \Z_{2j} \rtimes \{\pm 1\}$.
\end{thm}

\begin{proof}
It is routine to verify that the extension of each $(s,t) \in \Z_{2j} \rtimes\{\pm 1\}$ to the degree-$2$ vertices and each $\beta \in B$ are indeed automorphisms of $C_{2j}(i, j_\div p)$. Thus both $\Z_{2j} \rtimes \{\pm 1\}$ and $B$ are subgroups of $\aut(C_{2j}(i, j_\div p))$. 
It is similarly routine to verify that for all $a \in \Z_{2j}$ and $(s,t) \in \Z_{2j} \rtimes \{\pm 1\}$,
$ (s,t) \circ \beta_a = \beta_{s+ta} \circ (s,t)$,
so $B$ is a normal subgroup of $\aut(C_{2j}(i, j_\div p))$. Notice, $B \cap \Z_{2j} \rtimes \{\pm 1\} = \emptyset$. It now suffices to show that every automorphism of $C_{2j}(i, j_\div p)$ can be represented as the composition of an element of $B$ and an element of $\Z_{2j} \rtimes \{\pm 1\}$.

Let $\alpha \in \aut(C_{2j}(i, j_\div p))$. Then $\alpha$ restricts to a bijection $\alpha'$ on the set of degree-4 vertices, $\{u_a \mid a \in \Z_{2j}\}$; we will show that in fact $\alpha'\in \aut(C_{2j}(i,j))$. By definition $u_a$ and $u_b$ are adjacent in $C_{2j}(i,j)$ if and only if, as vertices in $C_{2j}(i, j_\div p)$, they are either adjacent or there are two vertex-disjoint paths of length $p+1$ between them, with all interior vertices of degree $2$. Such a relationship exists between $u_a$ and $u_b$ in $C_{2j}(i, j_\div p)$ if and only if it also exists between $\alpha(u_a)$ and $\alpha(u_b)$ in $C_{2j}(i, j_\div p)$, and this in turn holds if and only if $\alpha'(u_a)$ and $\alpha'(u_b)$ are adjacent in $C_{2j}(i,j)$. Hence $\alpha' \in \aut(C_{2j}(i,j))$.

If $C_{2j}(i,j)$ is twin-free, then $\alpha' \in \aut(C_{2j}(i,j)) = \Z_{2j} \rtimes H' = D_{4j}$ by
Theorem~\ref{thm:TwinFreeAutG}.
By Lemma~\ref{lem:whenTwins}, only two connected, two-generator circulant graphs satisfying $n=2j$ have twins, namely $C_4(1,2)$ and $C_6(1,3)$. As we argued in the proof of Theorem~\ref{thm:AutII}, in both of these cases, the edges corresponding to arcs of voltage $i=1$ induce a Hamilton ${2j}$-cycle in $C_{2j}(i,j)$. Moreover, $\alpha'$ respects this ${2j}$-cycle, because $\alpha$ must respect the Hamilton cycle in $C_{2j}(i, j_\div p)$ induced by all edges having no endvertex of degree $2$. Thus $\alpha'$ is an automorphism of an ${2j}$-cycle and so we can again conclude that $\alpha'\in D_{4j} =\Z_{2j} \rtimes H'$.

Thus in all cases, $\alpha' = (s, t)$ for some $s \in \Z_{2j}$ and $t \in \{\pm 1\}$. Since $\Z_{2j} \rtimes \{\pm 1\}$ is a subgroup of $\aut(C_{2j}(i, j_\div p))$, both $(s,t)$ and $(-s, t)$ are in $\aut(C_{2j}(i, j_\div p))$. Then $\gamma = \alpha \circ (-s, t)$ is an automorphism of $C_{2j}(i, j_\div p)$ that fixes every degree-$4$ vertex. We will show that $\gamma \in B$.

Let $a \in \Z_{2j}$. There are exactly two degree-$2$ vertices that are both adjacent to $u_a$ and distance $p+1$ from $u_{a+j}$, namely $v_a^1$ and $v_{a+j}^q$. Since $\gamma$ is a bijection that fixes both $u_a$ and $u_{a+j}$, $\gamma$ either fixes or interchanges $v_a^1$ and $v_{a+j}^q$.
For any $1 < r \le q$, $v_a^r$ is adjacent to $v_a^{r-1}$, so 
$\gamma(v_a^r)$ must be adjacent to $\gamma(v_a^{r-1})$. We can use this to prove inductively that for all such $r$,
\[
\gamma(v_a^r) = \begin{cases} v_a^r, \quad &\text{if }\gamma(v_a^1) = v_a^1,\\ 
v_{a+j}^{q+1-r}, & \text{if }\gamma(v_a^1) = v_{a+j}^q.
\end{cases}
\]
That is, $\gamma$ either fixes or interchanges the two paths between $u_a$ and $u_{a+j}$ that have all interior vertices of degree $2$. If it fixes them, let $e_a=0$, and if it switches them, let $e_a = 1$. Note that $e_a = e_{a+j}$. Then $
\gamma = \beta_0^{e_0} \beta_1^{e_1} \cdots \beta_{j-1}^{e_{j-1}} \in B.$
Since we defined $\gamma = \alpha\circ(-s, t)$, we have
$\alpha = \gamma \circ (s,t)$. This completes the proof that the automorphism group is the semidirect product $B \rtimes (\Z_{2j} \rtimes \{\pm 1\})$.
\end{proof}

\begin{prop}\label{prop:CaseIIIDet}
If $j \ge 3$, then a minimum determining set of $C_{2j}(i,j_\div p)$ is $W=\{v_0^1, v_1^1, \dots, v_{j-1}^1\}$ if $p \ge 2$ and 
$W= \{u_0, v_0^1, v_1^1, \dots, v_{j-1}^1\}$ if $p=1$.
\end{prop}

 \begin{proof} First note that for each $a \in \Z_n$, a determining set must contain at least one degree-2 vertex on the two paths between $u_a$ and $u_{a+j}$, for otherwise $\beta_a$ is a nontrivial automorphism fixing the set.

 Assume $p \ge 2$. The set $W=\{v_0^1, v_1^1, \dots, v_{j-1}^1\}$ contains exactly one degree-2 vertex on each pair of paths between $u_a$ and $u_{a+j}$. Thus any automorphism fixing this set cannot include any $\beta_a$, so it is of the form $\iota_B \circ(s, t) = (s, t)$. Since $u_0$ is the only degree-4 vertex adjacent to $v_0^1$, any $(s, t)$ fixing every vertex in $W$ also fixes $u_0$ and is therefore of the form $(0, t)$. By definition, $(0, t) \cdot (v_1^1) = v_t^1$, so the assumption that $(0,t)$ fixes every vertex in $W$ implies $t\equiv 1$. The automorphism $\iota_B \circ(0,1)$ is the identity of $\aut(C_{2j}(i,j_\div p)$ and so $W$ is determining.

 Next assume $p=1$. By Lemma~\ref{lem:whenTwinsIII}, $v_a^1$ and $v_{a+j}^1$ are twins for all $0 \le a < j$. Note that $\beta_a$ is an automorphism that interchanges these twin vertices and leaves all other vertices fixed. 
 However, no minimum twin cover is determining because it is fixed by the nontrivial automorphism $(\beta_0 \beta_1 \dots \beta_{j-1}) \circ (j, 1)$. On the other hand, the set $W=\{u_0, v_0^1, v_1^1, \dots, v_{j-1}^1\}$ is determining. Since $u_0 \in W$, any automorphism fixing $W$ is of the form $\beta \circ (0, t)$, where $\beta = \beta_0^{e_0}\cdot \beta_1^{e_1} \cdot \dots \cdot \beta_{j-1}^{e_{j-1}}$
for some $e_0, e_1, \dots , e_{j-1} \in \{0, 1\}$, and $t \equiv \pm 1$. If $t\equiv 1$, then $(0,1)$ fixes every element of $W$. If $e_a=1$ for any $a \in \{0, 1, \dots, j-1\}$, then $\beta\cdot(v_a^1) = v_{j+a}^1 \neq v_a^1$, so this would imply that $\beta = \iota_B$ and $\beta \circ(s,t)$ is the identity. So for this automorphism to be nontrivial, we require $t\equiv -1$. Then 
\[
 [\beta \circ (0, -1)] \cdot (v_1^1) = \beta \cdot(v_{2j-1}^1) =\begin{cases}
v_{2j-1}^1, \quad &\text{ if } e_1=0,\\
v_{j-1}^1, & \text{ if } e_1=1.
\end{cases}
\]
Since neither $2j-1$ nor $j-1$ equals $1$ (by the assumption that $j\ge 3$), there is no $\beta \in B$ for which $\beta\circ(0,-1)$ fixes every element of $W$, and and so $W$ is a determining set.
 \end{proof}

 \begin{thm}\label{thm:DistCostIII} 
 Assume $C_{2j}(i,j_\div p)$ is connected.
 \begin{enumerate}[(1)]
 \item If $p \ge 2$, then
$\det(C_{2j}(i,j_\div p)) =j, \, \dist(C_{2j}(i,j_\div p)) = 2,$ and
\[
\rho(C_{2j}(i,j_\div p))=\begin{cases}
 j+1,\quad &p = 2 \text{ and } j\in \{2, 3, 4, 5\}, \text{ or } p = j=3, \\
 j, &\text{otherwise.}
 \end{cases}\]
\item If $p = 1$, then $\det(C_4(1,2_\div 1)) = 4$ and $\dist(C_4(1,2_\div 1)) = 3$.
If $j \ge 3$, $\det(C_{2j}(i,j_\div 1)) = j+1$, $\dist(C_{2j}(i,j_\div 1)) = 2$, and
$\rho(C_{2j}(i,j_\div 1))=j+3.$
\end{enumerate}
\end{thm}

\begin{proof}
We first prove (2); assume $p=1$. 
For each $a \in \Z_{2j}$, $\beta_a$ interchanges twins $v_a^1$ and $v_{a+j}^1$, and fixes all other vertices. Thus $v_a^1$ and $v_{a+j}^1$ must be in different color classes in any distinguishing coloring.

 If $j=2$, we have $C_4(1,2_\div 1)$.
 The determining and distinguishing number are discussed in Example~\ref{ex:C12div1}. Note that the cost of 2-distinguishing is undefined in this case.

 Now suppose $j\ge 3$.
 By Proposition~\ref{prop:CaseIIIDet}, $\det(C_{2j}(i,p_\div 1)) = j+1$.
 From the proof of Theorem~\ref{thm:DetDistCostCaseI}, $u_i$ and $u_j$ are nonadjacent in $C_{2j}(i,j)$, so they are also nonadjacent in $C_{2j}(i,j_\div p)$. Color the vertices in $R = \{u_0, u_i, u_j, v_0^1, v_1^1\dots, v_{j-1}^1\}$ red and all other vertices blue, and assume $\beta \circ (s, t) \in B \rtimes (\Z_{2j} \rtimes \{\pm 1\})$ preserves these two color classes. 
 Among the degree-4 vertices in $R$, $u_0$ and $u_i$ are adjacent to each other but neither is adjacent to $u_j$. 
 Thus $[\beta \circ (s, t)]\cdot(u_j) = u_j$.
 This implies that $j \equiv s+tj \equiv s+j$, so $s \equiv 0$.
 Note that $[\beta\circ(0,t)]\cdot(u_0) = u_0$.
 Then $u_i$ must also be fixed by $\beta\circ(0,t)$
 so $ti \equiv i$ and so $t\equiv 1$. In order for $\beta\circ(0, 1)$ to fix all of the degree-2 vertices in $R$, $\beta = \iota_B$. Thus this is a $2$-distinguishing coloring with $|R|=j+3$ red vertices and (by the assumption that $j\ge 3$) at least $j+3$ blue vertices. Hence $\rho(C_{2j}(i,j_\div p)) \le j+3$.

 To find a lower bound on cost, first note that since $v_a^1$ and $v_{a+j}^1$ must be in different color classes for all $a \in \Z_{2j}$, each color class in a 2-distinguishing coloring contains at least $j$ vertices.
 Suppose there is just one degree-$4$ vertex, $u_a$, in the minimum size color class. We know that $(s, t) = (2a, -1)$ fixes $u_a$. For each $b\in \Z_{2j}$, $(2a, -1)$ exchanges $v_b^1$ and $v_{2a-b}^1$.
 If $v_b^1$ and $v_{2a-b}^1$ have opposite colors, then let $e_b=1$, and otherwise let $e_b=0$. By the first sentence of the paragraph, $e_{b+j} = e_b$. Then let 
 \[
 \beta = \prod \{\beta_b \mid e_b=1 \text{ and } 0 \le b <j\}.
 \]
 Then the nontrivial automorphism $\beta\circ(2a, -1)$ preserves the color classes, so the coloring is not distinguishing. Suppose instead that there are exactly two degree-4 vertices, $u_a$ and $u_b$, in the minimum color class. Using the same process as above, we can find $\beta \in B$ such that the nontrivial automorphism $\beta \circ (a+b,-1)$ preserves the color classes, so the coloring is not distinguishing. Thus $\rho(C_{2j}(i,j_\div p)) >j+2$. 
 
\smallskip

We now prove (1); assume $p \ge 2$. 
In this case $\det(C_{2j}(i, j_\div p)) = j$ by Proposition~\ref{prop:CaseIIIDet}.
To show that $C_{2j}(i,j_\div p)$ is always $2$-distinguishable, color the vertices in $R = \{u_0, v_0^1, v_1^1, \dots, v_{j-1}^1\}$ red and all other vertices blue and assume $\beta \circ (s,t)$ preserves the color classes. Since $u_0$ is the only degree-4 vertex in $R$, it must be fixed by $\beta \circ (s,t)$, so $s\equiv 0$. 
Recall that $H' = \{\pm 1\}$.
If $t\equiv -1$, 
\[
[\beta\circ(0, -1)] \cdot (v_1^1) = \beta \cdot (v_{2j-1}^1) = 
\begin{cases} 
v_{2j-1}^1, \, &\text{ if } e_{j-1} = 0,\\
v_{j-1}^p, &\text{ if } e_{j-1} = 1.
\end{cases}
\]
Since neither $v_{2j-1}^1$ nor $v_{j-1}^p$ is in $R$, no automorphism of the form $\beta\circ(0, -1)$ preserves $R$. 
Hence $t \equiv 1$. Since
 $(s,t) = (0, 1)$ fixes every degree-2 vertex in $R$, $\beta = \iota_B$, and so $\beta\circ(s, t)$ is the identity. This is therefore a $2$-distinguishing coloring, so $\dist(C_{2j}(i,j_\div p)) = 2$ and $\rho(C_{2j}(i,j_\div p)) \le j+1$.

To establish a lower bound on cost, note that in any $2$-distinguishing coloring, $v_a^r$ and $v_{a+j}^{p+1-r}$ must have different colors for at least one $1 \le r \le p$; otherwise $\beta_a$ preserves the color classes. This means that each color class must have at least $j$ vertices.

The remainder of the proof is devoted to establishing when we can find a $2$-distinguishing coloring with exactly $j$ vertices in a color class. To this end, we define a {\it set of representatives of $\Z_j$ in $\Z_{2j}$} to be any set of the form $S = \{a_0, a_1, \dots, a_{j-1}\}\subset \Z_{2j}$, where $a_k \in \{k, k+j\}$. 
In any $2$-distinguishing coloring, 
candidates for a color class of size $j$ in a $2$-distinguishing coloring are sets of the form
\begin{equation}\label{eqn:possibleMin}
R = \{v_{a_0}^{r_0}, v_{a_1}^{r_1}, \dots, v_{a_{j-1}}^{r_{j-1}}\},\end{equation}
where the subscripts 
$S = \{a_0, a_1, \dots, a_{j-1}\}$ form a set of representatives of $\Z_j$ in $\Z_{2j}$ and the superscripts $r_0, r_1, \dots , r_{j-1}$ are in $\{1, \dots, p \}$. 
If every such $R$ is preserved by some nontrivial automorphism, then the cost is $j+1$. Conversely, if we can find one such $R$ preserved only by the identity, then the cost is $j$.

Let $\beta \circ (s,t) \in \aut(C_{2j}(i,j_\div p)) = B \rtimes (\Z_{2j} \rtimes \{\pm 1\})$. 
Note that the action of $(s,t)$ on the elements of $R$ changes only the subscripts. More precisely, let $(s,t) \in \Z_{2j} \rtimes \{\pm 1\}$ act on the elements of $S$ as it would on elements of $\Z_{2j}$: $(s,t) \cdot (a_k) = s + ta_k$, which may or may not be an element of $S$.

Table~\ref{tab:PreservingAction} in Appendix~\ref{sec:BigTable} shows that for $j \in \{2, 3, 4, 5\}$, every set of representatives of $\Z_j$ in $\Z_{2j}$ has a nontrivial automorphism in $\Z_{2j} \rtimes \{\pm 1\}$ preserving it. 
If $p=2$, it can be shown that this implies that any set $R$ of the form in Equation~\ref{eqn:possibleMin} has a nontrivial automorphism preserving it; see Lemma~\ref{lem:preserveS} in Appendix~\ref{sec:LastDetails}. Thus for $p=2$ and $j\in \{2, 3, 4, 5 \}$, $\rho(C_{2j}(i,j_\div p)) = j+1$.

However, for $j \ge 6$, the set of representatives $S= \{0, 1, j+2, 3, 4, \dots, j-1\}$
is preserved only by the identity 
of $\Z_{2j} \rtimes \{\pm 1\}$; 
see Lemma~\ref{lem:CaseIIIj6} in Appendix~\ref{sec:LastDetails}.
Let $R=\{v_0^1, v_1^1, v_{j+2}^1, v_3^1, \dots, v_{j-1}^1\}$ and assume $\beta \circ (s, t)$ preserves $R$.
 For all $a \in S$, 
 \[
[\beta\circ(s, t)] \cdot (v_a^1) = \beta (v_{s+ta}^1) = 
\begin{cases} 
v_{s+ta}^1, \, &\text{ if } e_{s+ta} = 0,\\
v_{s+ta+j}^p, &\text{ if } e_{s+ta} = 1.
\end{cases}
\]
If $s+ta \in S$, then $v_{s+ta}^1 \in R$. However, by Lemma~\ref{lem:CaseIIIj6}, if $(s,t)$ is nontrivial, then for some $a \in S$, $s+ta\notin S$, so $v_{s+ta}^1 \notin R$. Since $S$ is a set of representatives, $s+ta+j \in S$, but $v_{s+ta+j}^p \notin R$. So for $R$ to be preserved, $(s, t) = (0, 1)$, which in turn implies $\beta = \iota_B$.
So for $p=2$ and $j \ge 6$, $\rho(C_{2j}(i,j_\div p)) = j$.

Now assume $p \ge 3$. Let $R=\{v_0^2, v_1^p\}$ if $j=2$ and $R = \{v_0^2, v_1^1, \dots, v_{j-2}^1, v_{j-1}^p\}$ if $j \ge 3$.
The subscripts 
form the set of representatives $S=\{0, 1, \dots, j-1\}$.
Assume $\beta \circ (s, t)$ preserves $R$. 
The only element of $R$ that is not adjacent to a degree-4 vertex is $v_0^2$, so it must be fixed by $\beta \circ (s, t)$, which implies that $s\equiv 0$ and $e_0=0$. Suppose that $t\equiv -1$. The image of $R$ under $(0, -1)$ is
\[
(0, -1)\cdot (R) = \{v_0^2, v_{2j-1}^1, \dots, v_{j+2}^1, v_{j+1}^p\} 
= \{v_0^2, v_{j+1}^p, v_{j+2}^1, \dots, v_{2j-1}^1 \}.
\]
The set of representatives of $\Z_j$ in $\Z_{2j}$ formed by the subscripts on these degree-$2$ vertices is $\{0, j+1, j+2, \dots, 2j-1\}$; to get the subscripts back to the set $S$, we must have $\beta = \beta_1 \beta_2 \dots \beta_{j-1}$. Then 
\[
[\beta\circ (0, -1)]\cdot (R) = \{v_0^2, v_1^1, v_2^p, \dots, v_{j-1}^p\}.
\]
If $j=2$, then $[\beta_1 \circ(0, -1)]\cdot (R) = \{v_0^2,v_1^1\} \neq R$. If $j>3$, then $v_2^p$ is in $[\beta\circ (0, -1)]\cdot (R)$ but not in $R$, contradicting the assumption that $\beta\circ(0, -1)$ preserves $R$. Thus for $j=2$ and $j>3$, we conclude $t\equiv 1$, which in turn implies $\beta = \iota_B$ and so the only automorphism preserving $R$ is the identity. In these cases, $\rho(C_{2j}(i,j_\div p))= j$.

If $j=3$, however, $[\beta_1 \beta_2 \circ(0, -1)] \cdot (R) = \{v_0^2,v_1^1, v_2^p\} = R$.
If $p\ge 4$, this problem is easily addressed by replacing $R$ with $R'=\{v_0^1, v_1^2, v_2^2\}$.
In $R'$, $v_0^1$ is the only vertex adjacent to a degree-$4$ vertex, and so it is fixed by any automorphism preserving $R'$. Note that
$
(0, -1)\cdot (R') = \{v_0^1, v_5^2, v_4^2\}.
$
To get the subscripts back to the set of representatives $\{0, 1, 2\}$, we must have $\beta = \beta_1 \beta_2$, but
$\beta_1 \beta_2 \circ (0, -1) (R') = \{v_0^1, v_2^{p-1}, v_1^{p-1}\} \neq R'.$
Hence $t\equiv1$, so $\beta = \iota_B$ and thus the only automorphism preserving $R'$ is the identity. So for $j=3$ and $p\ge 4$, $\rho(C_n(i,j_\div p))= 3 =j$. 
However, if $j = p = 3$, it can be shown that the cost of $2$-distinguishing cannot be $j$; see Lemma~\ref{lem:CaseIIIpj3} in Appendix~\ref{sec:LastDetails}.

\end{proof}

\section{Future Work}\label{sec:Open}

A natural extension of our work is to find symmetry parameters for connected $C_n(A)$ where $|A|>2$. 
Note that if $\overline A$ denotes the complement of $A$ in $\{1, 2, \dots, n\}$, then $C_n(\overline A) = \overline{C_n(A)}$. It is known that the determining number, distinguishing number and, if relevant, the cost of 2-distinguishing are equal for a graph and its complement. This means it suffices to find these symmetry parameters for connected $C_n(A)$, where $2< |A| \le n/2$.

For two-generator circulant graphs, we have found it fruitful to divide into cases depending on the presence of twins or co-twins. The obvious generalization of Theorem~\ref{thm:TwinFreeAutG} would be that if $C_n(A)$ is twin-free and co-twin-free, then $\aut (C_n(A))= \Z_n \rtimes H$ where $H = \aut(\mathbb Z_n, A)$. Equivalently, by Godsil's result, if $C_n(A)$ is twin-free and co-twin-free, then $\mathbb Z_n$ is a normal subgroup of $\aut(C_n(A))$.
For $C_n(A)$ with twins, we can use the results of Section~\ref{sec:TwinCase} to compute the symmetry parameters in terms of those of the twin quotient graph $\widetilde{C_n(A)}$.
By vertex transitivity, if one vertex in $C_n(A)$ has $k$ twins, then so does every vertex. 
It follows that the degree of every vertex is a multiple of $k$.
Together, these imply that the twin quotient graph is a twin-free circulant graph of order $n/k$ with fewer generators. A similar approach is helpful in considering circulant graphs with co-twins.

Another direction for future research would be to investigate the symmetry parameters of other subdivisions of connected $C_n(i,j)$, such as $C_n(i_\div \ell, j_\div m)$ when $i \neq j$ and when $i=j$ but $\ell \neq m$.

\section*{Acknowledgments}

The authors thank Debra Boutin, Puck Rombach, and especially Lauren Keough for valuable discussions in the early stages of this research. We also thank the anonymous reviewer for their extensive suggestions for improving this paper.

\bibliographystyle{abbrv}
\bibliography{ArXiV2}

\begin{appendices}

\section{Proposition~\ref{prop:autoFix0} under special conditions}\label{sec:AutoSpecialConditions}

 Recall that Proposition~\ref{prop:autoFix0} states that if $C_n(i,j)$ is twin-free and not $C_{10}(1,3)$, then any automorphism $\alpha$ fixing $0$ is an automorphism of the additive group $\Z_n$.
 The strategy of the proof is to use induction on $c+d = m$ to show that for all $0 \le c, d \in \Z$,
\begin{equation}\label{eqn:indApp} \alpha(ci+dj) = c\alpha(i) + d\alpha(j).\end{equation} 
The base case $m=1$ is trivial. 
In the body of the paper, we used an induction argument in the cases $j = n/2$ and $j < n/2$ with none of the special conditions in Lemma~\ref{lem:CommonNeighbors} holding. In this section, we assume $j < n/2$ and modify the induction argument to cover the cases where special conditions hold.
It can be verified computationally that the proposition holds in the two cases when two special conditions hold, namely, $C_{12}(3, 5)$ and $C_{12}(1, 3)$. In what follows, we assume just one special condition holds.

Assume $4i \equiv 0$. Applying $\alpha$ to $N(0) \cap N(2j) = \{j\}$ gives $N(0) \cap N(\alpha(2j)) = \{\alpha(j)\}$. 
From Lemma~\ref{lem:CommonNeighbors}, 
$0$ and $u$ have exactly one common neighbor only when $u = \pm 2j$, in which case the common neighbor is $\pm j$. 
 Since $\alpha$ is injective, 
$\alpha(j) \equiv \pm j,\, \alpha(-j) \equiv - \alpha(j), \, \alpha (2j) \equiv 2 \alpha(j).$
Since $N(0) = \{i, -i, j, -j\}$ and $\alpha$ fixes $0$, $N(0) = \{\alpha(i), \alpha(-i), \alpha(j), \alpha(-j)\}$. The last two elements of this set are $j$ and $-j$ (in some order), so $\alpha(i) \equiv \pm i$ and $\alpha(-i) \equiv -\alpha(i)$. 
Next, applying $\alpha$ to $N(0)\cap N(2i) = \{i, -i\}$ gives 
$N(0)\cap N(\alpha(2i)) = \{\alpha(i), \alpha(-i)\} = \{i, -i\}$. 
This means $\alpha(2i) \equiv\pm 2i$. We are assuming $4i\equiv 0$, so $2i \equiv -2i$ and thus $\alpha(2i) \equiv 2 \alpha(i)$.

Applying $\alpha$ to $N(0) \cap N(i+j) = \{i, j\}$ gives $N(0) \cap N(\alpha(i+j)) = \{\alpha(i), \alpha(j)\}$. From 
Lemma~\ref{lem:CommonNeighbors}, 
if $N(a)\cap N(w)= \{a+u, a+v\}$, where $u\in \{i, -i\}$ and $v\in \{j, -j\}$, then $w\equiv a+u+v$.
Thus $\alpha(i+j) \equiv \alpha(i) + \alpha(j)$. This completes establishing the result for the base case $m = c+d =2$.

For the inductive step, let $m = c+d \ge 3$ and assume Equation~\ref{eqn:indApp} holds for $m-1$ and $m-2$. Since $m \ge 3$ and $c, d \ge 0$, either $c\ge 2$ or $d\ge 2$. If $d \ge 2$, then from Lemma~\ref{lem:CommonNeighbors},
$
N(ci + (d-2)j) \cap N(ci+dj) = \{ci +(d-1)j\}.
$
Applying $\alpha$ and the inductive hypothesis gives
$
N(c\alpha(i) + (d-2)\alpha(j)) \cap N(\alpha(ci+dj)) = \{c\alpha(i) +(d-1)\alpha(j)\}.
$
When $4i\equiv 0$, if $N(u) \cap N(v) = \{w\}$, then $2w = u+v$, so
$
2\big (c\alpha(i) +(d-1)\alpha(j)\big ) \equiv c\alpha(i) + (d-2)\alpha(j)+\alpha(ci+dj),
$
which simplifies to $\alpha(ci+dj) \equiv c\alpha(i) + d\alpha(j)$.

Now assume $d < 2$. First suppose $d=0$. We showed above that $\alpha(i) = \pm i$ and $\alpha(-i) \equiv - \alpha(i)$. Under the condition $4i\equiv 0$, $-i\equiv 3i$. Substituting gives 
$
\alpha(3i) \equiv \alpha(-i) \equiv -\alpha(i) \equiv -(\pm i) \equiv -(\mp 3i) \equiv \pm 3i \equiv 3\alpha(i).
$
Thus $\alpha(ci) \equiv c\alpha(i)$ for all $c \in \{0, 1, 2, 3\}$, and these are all the values of $c$ that we need to verify under the condition $4i\equiv 0$. 

Last suppose $d=1$. From Lemma~\ref{lem:CommonNeighbors},
$
N((c-1)i) \cap N(ci+j) = \{ci, (c-1)i+j\}.
$
Applying $\alpha$ and the inductive hypothesis,
\begin{align*}
N((c-1)\alpha(i)) \cap N(\alpha(ci+j)) &= \{c\alpha(i), (c-1)\alpha(i)+\alpha(j)\} \\
& = \{(c-1)\alpha(i)+\alpha(i), (c-1)\alpha(i)+\alpha(j)\}.
\end{align*}
Letting $a = (c-1)\alpha(i)$, we can rewrite this as 
$
N(a) \cap N( \alpha(ci+j)) = \{a+ \alpha(i), a+\alpha(j)\}.
$
As noted earlier, in this situation Lemma~\ref{lem:CommonNeighbors} implies that 
\[
\alpha(ci+j) \equiv a + \alpha(i) + \alpha(j) \equiv (c-1)\alpha(i) + \alpha(i) + \alpha(j) \equiv c\alpha(i) + \alpha(j).
\]
This completes the proof when $4i\equiv 0$. By interchanging $i$ and $j$, we get a completely analogous proof under the special condition $4j\equiv 0$.

\smallskip

Next assume $3i\equiv j$.
For all $c\ge 0$, there exist $q\ge 0$ and $r \in \{0, 1, 2\}$ such that $c = r + 3q$. Thus
$ci = (r + 3q)i = ri + q(3i) = ri + qj$.
This implies that we only have to prove $\alpha(ri+qj) \equiv r \alpha(i) + q \alpha(j)$ for $r \in \{0, 1, 2\}$ and $q\ge 0$.
We modify our approach as follows: for each fixed $r\in\{0, 1, 2\}$, we use induction only on $q$.

Let $r=0$. The base case $q=1$ is trivial. For $q=2$, apply $\alpha$ to $N(0) \cap N(2j) = \{j\}$ to get $N(0) \cap N(\alpha(2j)) = \{\alpha(j)\}$. 
Similarly, $N(0) \cap N(\alpha(-2j)) = \{\alpha(-j)\}$. 
From Lemma~\ref{tab:CommonNeighbors}, there is only one situation in which two vertices have exactly one common neighbor. We can conclude
$\alpha(j) \equiv \pm j, \, \alpha(-j) \equiv -\alpha(j), \, \alpha(2j) \equiv 2 \alpha(j),\, \alpha(-2j) \equiv -2\alpha(j)$.
Now let $q \ge 3$ and assume the result holds for $q-1$ and $q-2$. 
Then applying $\alpha$ to $N((q-2)j) \cap N(qj) = \{(q-1)j\}$ gives $N((q-2)\alpha(j)) \cap N(\alpha(qj)) = \{(q-1)\alpha(j)\}$. Because the intersection is a singleton, Lemma~\ref{tab:CommonNeighbors} gives us $(q-2)\alpha(j) + \alpha(qj) = 2(q-1)\alpha(j)$, which simplifies to $\alpha(qj) = q\alpha(j)$, as desired.

Next let $r=1$. The first base case is now $q=0$, but it is still trivial. Since $\{\alpha(j), \alpha(-j)\} = \{j, -j\}$, $\alpha(i) \equiv \pm i$ and $\alpha(-i) \equiv - \alpha(i)$. Applying $\alpha$ to $N(0) \cap N(i+j) = \{i, j\}$ gives $N(0) \cap N(\alpha(i+j)) = \{\alpha(i), \alpha(j)\}$. When $3i\equiv j$, $|N(0) \cap N(v)| = 2$ 
only if either $v= i+j$, in which case $N(0) \cap N(v) = \{i,j\}$, or if $v=-i-j$, in which case $N(0) \cap N(v) = \{-i,-j\}$.
We conclude both that 
$\alpha(i+j) \equiv \alpha(i) + \alpha(j)$ and that $\alpha(i) \equiv i$ if and only if $\alpha(j) \equiv j$.
Thus $\alpha(i+qj) \equiv \alpha(i) + q\alpha(j)$ for the second base case $q=1$.

Let $q\ge 2$ and assume $\alpha(i+(q-1)j) \equiv \alpha(i) + (q-1)\alpha(j)$. From Lemma~\ref{lem:CommonNeighbors},
$
N((q-1)j)\cap N(i+qj) = \{i+(q-1)j, qj\}.
$
Applying $\alpha$, the inductive hypothesis and the result for $r=0$,
$N((q-1)\alpha(j))\cap N(\alpha(i+qj)) = \{\alpha(i)+(q-1)\alpha(j), q\alpha(j)\}.$
If $\alpha(i) \equiv i$ and $\alpha(j)\equiv j$, then substituting in yields
$
N((q-1)j)\cap N(\alpha(i+qj)) = \{i+(q-1)j, qj\}= \{i+(q-1)j, j+ (q-1)j\}.
$
From Lemma~\ref{lem:CommonNeighbors}, $\alpha(i+qj) \equiv i + qj \equiv \alpha(i) + q\alpha(j)$. If $\alpha(i)\equiv-i$ and $\alpha(j) \equiv -j$, 
$
N(-(q-1)j)\cap N(\alpha(i+qj)) = \{-i-(q-1)j, -qj\}= \{-i-(q-1)j, -j- (q-1)j\}.
$
Lemma~\ref{lem:CommonNeighbors} implies $\alpha(i+qj) \equiv -i -(q-1)j \equiv \alpha(i) + q \alpha(j)$.

Finally, let $r=2$. 
For the first base case $q=0$, apply $\alpha$ to $N(0)\cap N(2i)$ to get
$
N(0) \cap N(\alpha(2i)) = \{\alpha(i), \alpha(-i),\, \alpha(j)\} = \{i, -i, \alpha(j)\}.$
If $\alpha(i) \equiv i$ and $\alpha(j) \equiv j$, then from Lemma~\ref{lem:CommonNeighbors}, $\alpha(2i) \equiv 2i \equiv 2\alpha(i)$. 
If $\alpha(i) \equiv -i$ and $\alpha(j) \equiv -j$, then by Lemma~\ref{lem:CommonNeighbors}, $\alpha(2i) \equiv i-j$. Under the assumption $3i\equiv j$, this also gives $\alpha(2i) \equiv -2i \equiv 2\alpha(i)$.

Let $q\ge 1$ and assume $\alpha(2i+(q-1)j) \equiv 2\alpha(i) + (q-1)\alpha(j)$. Apply $\alpha$, the inductive hypothesis and the result for $r=1$ to 
$
N(i+(q-1)j) \cap N(2i+qj) = \{2i +(q-1)j, i+qj\}
$; we obtain
\[
N(\alpha(i) + (q-1)\alpha(j)) \cap N(\alpha(2i+qj)) = \{ 2\alpha(i) +(q-1)\alpha(j), \alpha(i)+q \alpha(j))\}.
\]
We can again examine the two possible cases 
$\alpha(i) \equiv i$, $\alpha(j) \equiv j$ and $\alpha(i)\equiv-i$, $\alpha(j) \equiv -j$ separately to conclude $\alpha(2i+qj) \equiv 2 \alpha(i) + q\alpha(j)$ for all $q \ge 0$. 
This concludes the proof under the special condition $3i\equiv j$. 

The proofs under the remaining special cases ($3i\equiv -j$, $3j\equiv i$ and $3j \equiv -i$) are just minor variations of the case $3i \equiv j$. We note that for the special cases $3i \equiv -j$ and $3j \equiv -i$, we need to prove the statement for $q \in \mathbb{Z}$ so the induction is on $|q|$. For $r=1$, this requires the second base cases of $q = \pm 1$; all other modifications are straightforward.

\section{Preserving sets of representatives of $\Z_j$ in $\Z_{2j}$}\label{sec:BigTable}

Table~\ref{tab:PreservingAction} shows that for $j \in \{2, 3, 4, 5\}$, every set $S$ of representatives of $\Z_j$ in $\Z_{2j}$ is preserved by a nontrivial automorphism of the form $(s, -1)\in \Z_{2j}\rtimes \{\pm 1 \}$.
\small
\begin{table}[h]
\begin{center}
\renewcommand{\arraystretch}{1.25}
\begin{tabular}{|c|c|c|}
\hline
 $j$ & $(s, -1)$ & $(a_0, \dots, a_{j-1})$ preserved by $(s, -1)$\\ 
 \hline 
 2 & $(1, -1)$ & $(0, 1), \, (2,3)$ \\
 \cline{2-3}
 & $(3, -1)$ & $(2, 1), \, (0, 3)$ \\
 \hline
 3 & $(0, -1)$ & $(0, 1, 5),\, (3,4,2)$ \\ 
 \cline{2-3} 
 & $(2, -1)$ & $(0, 1, 2), \, (3, 1, 5),\, (3, 4, 5)$\\
 \cline{2-3}
 & $(4, -1)$ & $(3, 1, 2), \, (0, 4, 2), \, (0, 4, 5)$\\
 \hline 
 4 &$(1, -1)$ & $(0, 1, 6, 3),\, (0, 1, 2, 7),\, (4, 5, 6, 3), \, (4, 5, 2, 7)$ \\ \cline{2-3} 
 & $(3, -1)$ & $(0, 1, 2, 3), \, (4, 1, 2, 7),\, (0, 5, 6, 3), \, (4, 5, 6, 7)$\\
 \cline{2-3} 
 & $(5, -1)$ & $(4, 1, 2, 3), \, (0, 5, 2, 3), \, (4, 1, 6, 7),\, (0, 5, 6, 7)$\\
 \cline{2-3} 
 & $(7, -1)$ & $ (4, 5, 2, 3), \, (4, 1, 6, 3),\, (0, 5, 2, 7),\, (0, 1, 6, 7)$\\
 \hline 
 5 &$(0, -1)$ & $(0, 6, 7, 3, 4),\, (0, 6, 2, 8, 4),\, (0, 1, 7, 3, 9), \, (0, 1, 2, 8, 9)$, \\
 & & $(5, 6, 7, 3, 4), \, (5, 6, 2, 8, 4), \ (5, 1, 2, 8, 9)$ \\ 
 \cline{2-3} 
 & $(2, -1)$ & $(0, 1, 2, 8, 4), \, (0, 1, 2, 3, 9),\, (0, 6, 2, 3, 9), \, (5, 1, 7, 8, 4)$, \\
 & & $(5, 1, 7, 3, 9), \, (5, 6, 7, 8, 4), \, (5, 6, 7, 3, 9)$\\
 \cline{2-3}
 & $(4, -1)$ & $(0, 1, 2, 3, 4), \, (0, 1, 7, 3, 4), \, (5, 1, 2, 3 , 9)$,\\
 & & $(0, 6, 7, 8, 4), \, (5, 6, 2, 8, 9), \, (5, 6, 7, 8, 9)$\\
 \cline{2-3}
 & $(6, -1)$ & $(5, 1, 2, 3, 4), \, (0, 6, 2, 3, 4),\, (5, 1, 2, 8, 4)$,\\ 
 & & $(0, 6, 7, 3, 9), (5, 1, 7, 8, 9), \, (0, 6, 7, 8, 9)$\\
 \cline{2-3}
 & $(8, -1)$ & $(5, 6, 2, 3, 4), \, (5, 1, 7, 3, 4),\, (0, 1, 7, 8, 4)$,\\
 & & $(5, 6, 2, 3, 9), \,(0, 6, 2, 8, 9),\, (0, 1, 7, 8, 9)$\\
 \hline
\end{tabular}
\caption{$(s, -1)$ preserving $(a_0, \dots, a_{j-1})$ .}\label{tab:PreservingAction}
\end{center}
\end{table}
\normalsize

\section{
Lemmas for Theorem~\ref{thm:DistCostIII}
}\label{sec:LastDetails}
\begin{lemma}\label{lem:preserveS}
Let $p=2$. Suppose that each set of representatives of $\Z_j$ in $\Z_{2j}$ is preserved by some nontrivial automorphism in $\Z_{2j} \rtimes \{\pm 1\}$. Let $R = \{v_{a_0}^{r_0}, v_{a_1}^{r_1}, \dots, v_{a_{j-1}}^{r_{j-1}}\}$, where 
$S = \{a_0, a_1, \dots, a_{j-1}\}$ is a set of representatives of $\Z_j$ in $\Z_{2j}$ and $r_0, r_1, \dots , r_{j-1} \in \{1, 2\}$. Then $R$ is preserved by a nontrivial automorphism in $\aut(C_{2j}(i,j_\div p)) = B \rtimes (\Z_{2j} \rtimes \{\pm 1\})$.
\end{lemma}

\begin{proof}
If $r_0 = r_1 = \dots =r_{j-1}$, then any $(s,t)$ preserving $S$ will preserve $R$. If the superscripts are not all equal, then let $\beta$ be the product of the $\beta_a$ corresponding to elements of $R$ with superscripts $r=2$; equivalently,
$\beta= \prod_{k=1}^{j-1} \beta_k^{r_k-1}$.
Applying $\beta$ to each element of $R$ creates a set
$\beta(R)$ in which every superscript is $1$ and the subscripts form a new set of representatives of $\Z_j$ in $\Z_{2j}$, which we name $\beta(S)$. By assumption, there exists a nontrivial automorphism $(s, t)$ that preserves $\beta(S)$ and thus preserves $\beta (R)$. Applying $\beta$ to this set takes it back to $R$. Thus $\beta \circ (s, t) \circ \beta$ preserves $R$. Using $(s,t) \circ \beta_a = \beta_{s+ta} \circ (s,t)$, we can rewrite this as $\beta'\circ (s,t) \in B \rtimes (\Z_{2j} \rtimes \{\pm 1\})$, for some $\beta' \in B$.
\end{proof}

\begin{Ex}
Let $j=4$ and let $R= \{v_0^2, v_1^1, v_2^2, v_3^2\}$. The first, third and fourth vertices have superscripts of $2$, so let $\beta = \beta_0 \beta_2 \beta_3$. Then 
$
\beta(R) = \{v_4^1, v_1^1, v_6^1, v_7^1\}.
$
From Table~\ref{tab:PreservingAction}, $(5, -1)$ preserves the set of representatives $\beta(S) =\{4, 1, 6, 7\}$, so
$
 \left [(5, -1) \circ \beta\right ](R) = \{v_4^1, v_1^1, v_6^1, v_7^1\}.
$
Applying $\beta = \beta_0 \beta_2 \beta_3$ to each element of this set takes us back to $R$. Note that 
\[ (5, -1) \circ \beta = (5, -1)\circ \beta_0 \beta_2 \beta_3 
= \beta_5 \beta_3 \beta_2\circ (5, -1) = \beta_1 \beta_2 \beta_3 \circ (5, -1).\] Hence a nontrivial automorphism preserving $R$ is 
\[
\beta \circ (5, -1) \circ \beta = (\beta_0 \beta_2 \beta_3) \circ [\beta_1 \beta_2 \beta_3 \circ (5, -1)] = \beta_0 \beta_1 \circ(5, -1).
\]
\end{Ex}

\begin{lemma}\label{lem:CaseIIIj6} 
For all $j \ge 6$, the only $(s, t)\in \Z_{2j} \rtimes \{\pm 1\}$ preserving the set of representatives $S= \{0, 1, j+2, 3, 4, \dots, j-1\}$ is the identity. 
\end{lemma}

\begin{proof}
It suffices to show that for each $(s, t) \neq (0, 1)$, there is some $a \in S$ such that $(s,t) \cdot a \in \overline S = \{j, j+1, 2, j+3, j+4, \dots, 2j-2, 2j-1\}$. 
Since $j \ge 6$, $j+4 < 2j - 1$.

If $s \in \overline S$, then 
$(s,t) \cdot 0 = s \in \overline S$ and we are done. So assume $s \in S$.
The first three possible values of $s$ are handled below:
\begin{align*}
s=0:\quad & (0, -1) \cdot 1 = \, 2j-1,\\
s=1:\quad & (1, -1) \cdot 3 = 2j-2 \text{ and } 
(1, 1) \cdot 1 = 2,\\
s = j+2: \quad & (j+2, -1)\cdot 1 = j+1 \text{ and }
(j+2, 1) \cdot 1 = j+3.
\end{align*}
For the remaining values of $s$ and $t= -1$, note that $(s,-1)\cdot(j-1) = s-j+1\equiv s+j+1$. If $3 \le s \le j-2$ then $(s, -1) \cdot (j-1) \in \{j+4, \dots, 2j-1\} \subset \overline S$.
For $s = j-1$, $(j-1, -1) \cdot (j-3) = 2\in \overline S$. Similarly, for $t =1$, if $4 \le s \le j$, then 
$(s,1) \cdot (j-1) = s+j-1\in \{j+3, \dots 2j-1 \} \subset \overline S$.
For $s=3$, note that $(3, 1) \cdot (j+2) = 3 + (j+2) = j+4 \in \overline S$.
\end{proof}

\begin{lemma}\label{lem:CaseIIIpj3}
Suppose $j=p=3$. Let $R = \{v_{a_0}^{r_0}, v_{a_1}^{r_1}, v_{a_2}^{r_2}\}$, where $S= \{a_0, a_1, a_2\}$ is a set of representatives of $\Z_3 $ in $\Z_6$ and $r_0,r_1, r_2 \in \{1, 2, 3\}$. Then there is a nontrivial automorphism of $C_{2j}(i, j_\div p) = C_6(i, 3_\div 3)$ that preserves $R$.
\end{lemma}

\begin{proof} 
If $r_0= r_1 = r_2$, then $R$ is preserved by the nontrivial automorphism preserving $S$ in Table~\ref{tab:PreservingAction}. 
Assume two subscripts are equal. 
If $r_0,r_1, r_2 \in \{1, 3\}$, 
then we can use the same technique as in the proof of Lemma~\ref{lem:preserveS} for the case $p=2$ when superscripts not all equal.

However, a new technique is required if $r_0,r_1, r_2 \in \{1, 2\}$ 
or $r_0,r_1, r_2 \in \{2, 3\}$, because $\beta_{a} (v_{a}^2) = v_{a+j}^2$; that is, 
applying $\beta_{a}$ does not change the superscript.
Let $R_1 = \{v_{a_0}^r, v_{a_1}^r, v_{a_2}^2\}$ and $R_2 = \{v_{a_0}^2, v_{a_1}^2, v_{a_2}^r\}$ where $r \in \{1, 3\}$.
The nontrivial automorphism $(a_0+a_1, -1)$ interchanges the subscripts $a_0$ and $a_1$, and $(a_0+a_1, -1) \cdot(a_2) = a_0+a_1 -a_2$.
Since $S$ is a set of representatives of $\Z_3$ in $\Z_6$, in $\Z_3$ we have 
$a_0+a_1 - a_2 \equiv a_0+a_1+2a_2 \equiv (a_0+a_1+a_2) + a_2 \equiv a_2$.
Thus in $\Z_6$, $a_0+a_1 - a_2 \equiv a_2$ or $3+ a_2$.
If $a_0+a_1-a_2 \equiv a_2$ in $\Z_6$, then both $R_1$ and $R_2$ are preserved by $(a_0+a_1, -1)$. If $a_0+a_1-a_2 \equiv 3+ a_2$ in $\Z_6$, then $R_1$ is preserved by $\beta_{a_2}\circ (a_0+a_1, -1)$ and $R_2$ is preserved by $\beta_{a_0} \beta_{a_1} \circ (3+ a_0+a_1, -1)$. 

Finally, we consider the case in which the three superscripts are different: without loss of generality, let $R_3 = \{v_{a_0}^1, v_{a_1}^3, v_{a_2}^2\}$.
If $a_0+a_1-a_2 \equiv a_2$, then $\beta_{a_0} \beta_{a_1} \beta_{a_2} \circ (3+ a_0+a_1, -1)$ preserves $R_3$. If $a_0+a_1-a_2 \equiv 3+ a_2$, then $\beta_{a_0} \beta_{a_1} \circ (3+a_0+a_1,-1)$ preserves $R_3$.
\end{proof}

\end{appendices}

\end{document}